\documentclass{article}

\usepackage[english]{babel}

\usepackage[letterpaper,top=2cm,bottom=2cm,left=3cm,right=3cm,marginparwidth=1.75cm]{geometry}
\usepackage[colorinlistoftodos,prependcaption]{todonotes}
\usepackage{xargs} 
\usepackage{amsmath}        
\usepackage{amssymb}
\newcommandx{\farhad}[2][1=]{\todo[inline,linecolor=red,backgroundcolor=red!25,bordercolor=red,#1]{Farhad: #2}}
\usepackage{graphicx}
\usepackage{amsfonts}
\usepackage{mathtools}

\usepackage{color}
\usepackage[normalem]{ulem}
\usepackage{amsmath}
\usepackage{algorithm}
\usepackage{algpseudocode}
\usepackage{subfig}

\usepackage[compress]{cite}
\usepackage{bm}
\usepackage{amsthm}
\usepackage{commath}
\def\*#1{\mathbf{#1}}
\newtheorem{thm}{Theorem}[section]

\newtheorem{lemma}[thm]{Lemma}

\newtheorem{remark}{Remark}
\newtheorem{defi}{Definition}

\newtheorem{assumption}{Assumption}
\newtheorem{prop}{Proposition}
\DeclareMathOperator{\nullop}{null}
\DeclareMathOperator{\spanop}{span}

\usepackage{amsmath}
\usepackage{graphicx}
\usepackage[colorlinks=true, allcolors=blue]{hyperref}

\title{\LARGE \bf
Two-timescale EXTRA for \\ Distributed Smooth Non-convex Optimization}

\author{Zeyu Peng, Farhad Farokhi, and Ye Pu
\thanks{The work of Z.~Peng is supported by the Australian Government Research Training Program Scholarship and the Rowden White Scholarship. The work of F.~Farokhi is supported by DP210102454 from the Australian Research Council (ARC). The work of Y.~Pu is supported by DE220101527 from the Australian Research Council (ARC).}
\thanks{The authors are with the Department of Electrical and Electronic Engineering at the University of Melbourne, Australia.}
}

\begin{document}
\maketitle

\begin{abstract}
In this paper, we study distributed optimization with smooth non-convex local objectives. We propose a novel variant of the well-known EXact firsT-ordeR Algorithm (EXTRA), called Two-timescale EXTRA, by introducing two distinct step-sizes. Leveraging the two-timescale strategy, we construct a Lyapunov function and establish the sub-linear convergence of Two-timescale EXTRA to a consensual first-order stationary point. Additionally, we introduce an off-line sequential method for algorithm parameter selection, and the numerical results support the theoretical guarantees.
\end{abstract}

\section{INTRODUCTION}

In this paper, we consider the optimization problem:
\begin{align} \label{P1}
        \min_{\mathbf{x} \in \mathbb{R}^{p}} \left\{F(\mathbf{x}) = \sum_{i=1}^{n} f_{i}(\mathbf{x})\right\}
\end{align}
where $\mathbf{x}$ represents the decision vector and each $f_{i}$ is differentiable and smooth but potentially non-convex. The goal is to utilize $n$ agents to iteratively and collaboratively solve the optimization problem \eqref{P1} in a distributed manner. We assume these $n$ agents are connected by an undirected and connected network $\mathcal{G} = \{\mathcal{V},\mathcal{E} \}$. Each agent $i \in \mathcal{V} = \{1,2,\dots,n\}$ only has information about its local objective function $f_{i}$ and local estimate $\mathbf{x}_{i}$. During each iteration, agents can exchange information with their neighbors within the network, allowing all agents to collectively work towards solving the target optimization problem \eqref{P1}. The distributed optimization problem over a connected network is widely exists in engineering applications such as learning, estimation and control~\cite{facchinei2015parallel,hong2017prox,gao2019reinforcement,gan2012optimal,ram2009distributed}.

Distributed optimization algorithms for solving the optimization problem in~(\ref{P1}) can be categorized into two main classes: dual-decomposition methods and consensus-based methods. Dual decomposition approaches aim to minimize an augmented Lagrangian that incorporates agreement-enforcing constraints among agents, typically through iterative updates of both primal and dual variables. When $f_{i}$ is (strongly) convex, the existing algorithms include Decentralized Alternating Directions Method of Multipliers (DADMM) \cite{chang2014multi,schizas2007consensus,shi2014linear}, Decentralized Linearized alternating direction Method of multipliers (DLM)~\cite{ling2015dlm} and Incremental Primal-Dual (PD) Gradient method~\cite{alghunaim2020linear} with linear convergence rate were established in~\cite{alghunaim2020linear,shi2014linear}. For non-convex settings, the theoretical convergence guarantee of Alternating Directions Method of Multipliers (ADMM) was established in~\cite{wang2019global}. Several variants of ADMM have been proposed to improve its applicability and convergence properties, including the Proximal Primal-Dual Algorithm (PROX-PDA) \cite{hong2016decomposing}, Linearized ADMM~\cite{lu2021linearized}, Gradient Primal-Dual Algorithm (GPDA)~\cite{hong2018gradient} and DecentrAlized Primal-Dual (ADAPD)\cite{10035454}. Among these, sub-linear convergence to an exact minimizer has been established for PROX-PDA~\cite{hong2016decomposing} and ADAPD~\cite{10035454}, whereas explicit convergence rates for the remaining algorithms have not been provided.

Consensus-based methods involve a consensus step to average local estimates across all agents \cite{xiao2004fast}, combined with a gradient descent step to guarantee a decreasing objective function at each iteration. Distributed Gradient Descent (DGD)~\cite{nedic2009distributed} is among the earliest distributed optimization methods, which requires a diminishing step-size to achieve exact convergence in both convex~\cite{yuan2016convergence} and non-convex settings~\cite{zeng2018nonconvex}.  The EXact firsT-ordeR Algorithm (EXTRA)~\cite{shi2015extra}, Distributed Inexact Gradient trackING (DIGING)~\cite{nedic2017achieving} and Near-DGD~\cite{iakovidou2021convergence} are three consensus based algorithms to address the convergence limitation of DGD based on three different techniques. EXTRA~\cite{shi2015extra} employs two distinct mixing matrices, which generate a non-zero correction term to eliminate the need for a vanishing step-size. When each local objective function $f_{i}$ is (strongly) convex, EXTRA achieves (linear) sub-linear convergence to the exact minimizer. DIGING~\cite{nedic2017achieving} is a gradient tracking-based method that introduces an auxiliary local variable to estimate the global gradient. However, it requires the exchange of two local variables at each iteration, resulting in increased communication overhead. The linear and sub-linear convergence to an exact minimizer were established for DIGING in convex~\cite{nedic2017achieving} and non-convex~\cite{daneshmand2020second} settings. Near-DGD enhances DGD by incorporating an inner consensus loop—performing multiple communications at each iteration—to achieve better agreement among agents. The convergence to an exact minimizer has been established in non-convex setting \cite{iakovidou2021convergence}. Network succEssive conveX approximaTion (NEXT) algorithm~\cite{di2015distributed} is a gradient tracking-based distributed algorithm designed to solve problems involving a smooth objective function combined with a non-smooth but convex regularizer. Its extension, Inexact NEXT (In-NEXT)~\cite{di2016distributed}, generalizes the framework to accommodate time-varying communication networks. However, similar to DIGING, these two algorithms both require to exchange two local variables at each iteration. Among these methods, EXTRA~\cite{shi2015extra} offers a favorable balance between theoretical convergence guarantees and practical implementation efficiency (with simple iterations and low communication cost). By incorporating an additional mixing matrix, it effectively addresses the convergence gap of DGD in the convex case without incurring extra communication overhead. Furthermore, the use of two separate mixing matrices makes EXTRA a highly flexible framework. Several algorithms such as DIGING~\cite{nedic2017achieving}, Distributed Proportional-Integral (PI)~\cite{yao2018distributed} and GPDA~\cite{hong2018gradient} are all closely related to EXTRA by defining a specific relationship between the two mixing matrices used in EXTRA. We also want to highlight that EXTRA strikes a balance between dual-decomposition-based and consensus-based methods, as its correction step can be interpreted as a form of dual variable update. In the convex case, some primal dual analysis of EXTRA has been established by \cite{mokhtari2016dqm,mokhtari2015decentralized}. In the non-convex setting, a unified framework that covers EXTRA was proposed in \cite{alghunaim2022unified}. However, their convergence analysis relies on the assumption that the two mixing matrices are linearly dependent.

In this work, we propose Two-timescale EXTRA (TT-EXTRA), a novel algorithm designed to solve distributed smooth non-convex optimization problems. TT-EXTRA can be viewed as a variant of EXTRA, distinguished by the incorporation of an additional step-size parameter. This design enables the construction of a Lyapunov function, through which we show that TT-EXTRA converges to the set of first-order stationary points of problems~\eqref{P1}. Compared to the convergence established in \cite{alghunaim2022unified}, our analysis is conducted under weaker assumptions (see Assumption~\ref{AssM}), which allows greater flexibility in the selection of mixing matrices. Based on the Lyapunov function, we prove that TT-EXTRA achieves a sublinear convergence rate that matches the complexity lower bound established in \cite{sun2019distributed}. Furthermore, we introduce an off-line sequential method for parameter selection, including the two step-sizes and mixing matrices. Finally, our numerical experiments demonstrate that the two-timescale scheme, together with the flexibility in selecting mixing matrices, yields improved practical performance. 
\subsection{Notations}
Throughout this paper, we assume that each agent $i$ maintains a local estimate of the global variable $\mathbf{x}$, denoted by $\mathbf{x}_i \in \mathbb{R}^p$, where $\mathbf{x}_i^r$ indicates its value at iteration $r$. For a given communication graph $\mathcal{G}$, $\mathcal{N}_i$ denotes the set of all neighbors of agent $i$ (i.e., the set of agents $j$ capable of transmitting their local estimates $\mathbf{x}_j^r$ to agent $i$). Matrices are represented by uppercase letters; for a matrix $A$, $\|A\|$ refers to its spectral norm, $\lambda_{\max}(A)$ to its largest eigenvalue, and $\lambda_2(A)$ to its second largest eigenvalue. Bold lowercase letters represent vectors; for a vector $\mathbf{x}$, $\|\mathbf{x}\|$ denotes the Euclidean norm, while $\|\mathbf{x}\|_A^2 := \mathbf{x}^T A \mathbf{x}$. Lastly, $\bm{1}$ denotes the all-ones vector and $\otimes$ denotes the Kronecker product.


\section{Two Timescale EXTRA (TT-EXTRA)} \label{sec:TT-EXTRA} 
In this section, we introduce the TT-EXTRA algorithm and its convergence properties. As shown in~\cite{10035454,shi2015extra}, if the communication network $\mathcal{G}$ is connected, the optimization problem~\eqref{P1} is equivalent to:
\begin{subequations}\label{P2}
    \begin{align}
    \min_{\mathbf{x} \in \mathbb{R}^{np}}& \left\{f(\mathbf{x}) = \sum_{i=1}^{n} f_{i}(\mathbf{x}_{i}) \right\},\\
    \mathrm{s.t.}\;& (W\otimes I_{p})\mathbf{x}=\mathbf{x},  
\end{align}
\end{subequations}
where $\mathbf{x} = (\mathbf{x}_{1}^T,\mathbf{x}_{2}^T,\dots, \mathbf{x}_{n}^T)^T \in \mathbb{R}^{np}$ and mixing matrix $W \in \mathbb{R}^{n\times n}$ satisfies $\nullop \{I_{n} \!-\! W \} = \spanop\{\bm{1}\}$. Several common choices of $W$ are presented in \cite{shi2015extra}.

\begin{defi} \label{def}
 A vector $\mathbf{x}^* \in$ $\mathbb{R}^{np}$ is said to be a first-order consensual stationary point of problem \eqref{P2} if it satisfies:
  \begin{itemize}
      \item $((I_n-W)\otimes I_p)\mathbf{x}^* = 0$ (consensus);
      \item $\bm 1^T \nabla f(\mathbf{x}^*) = 0$ (stationarity).
  \end{itemize}
  The set of all first-order consensual stationary points is denoted by $\mathcal{X}^*$.
\end{defi}

Given that $\nullop(I_{n} - W) = \spanop \{\bm 1\}$, the consensus condition in Definition~\ref{def} indicates that  $\mathbf{x}_{i}^*=\mathbf{x}_{j}^*$, $\forall i,j$, in $\mathbf{x}^*= (\mathbf{x}_{1}^*,\mathbf{x}_{2}^*,\dots, \mathbf{x}_{n}^*)^T$. Combined with the optimality condition, this implies that each element of $\mathbf{x}^*$ is a first-order stationary point of \eqref{P1}. We later show that the TT-EXTRA algorithm converges asymptotically to the set of first-order consensual stationary points $\mathcal{X}^*$.

\begin{algorithm}
\caption{TT-EXTRA}\label{alg:TTEXTRA}
\begin{algorithmic}[1]
    \State Choose $\beta, \rho > 0$ and mixing matrices $W, \Tilde{W} \in \mathbb{R}^{n\times n}$
    \State Choose $\mathbf{x}_{i}^{0}\in \mathbb{R}^{p}$ randomly, $\forall i\in\{1,\dots,n\}$
    \State Exchange $\mathbf{x}_{i}^{0}$ with agent $j \in \mathcal{N}_{i}$
    \State 	Let $\mathbf{y}_{i}^{0}=\rho \sum_{j=1}^{n} (\Tilde{W}_{ij}-W_{ij}) \mathbf{x}_{j}^{0}$
    \For{$r=0,1,2,\dots$}
        \For{$i=1,2,\dots,n$}
            \State $\mathbf{x}_{i}^{r+1} \!=\! (1\!-\!\frac{\rho}{\beta}) \mathbf{x}_{i}^{r}\!-\!\frac{1}{\beta} \nabla f_{i}(\mathbf{x}_{i}^{r}) \!+\! \frac{\rho}{\beta} \sum_{j\!=\!1}^{n} \!\Tilde{W}_{ij}\mathbf{x}_{j}^{r} \!-\! \frac{1}{\beta} \mathbf{y}_{i}^{r}$
            \State Exchange $\mathbf{x}_{i}^{r+1}$ with agent $j \in \mathcal{N}_{i}$
            \State $\mathbf{y}_{i}^{r+1} = \mathbf{y}_{i}^{r} + \rho \sum_{j=1}^{n} (\Tilde{W}_{ij}-W_{ij}) \mathbf{x}_{j}^{r+1} $
        \EndFor
    \EndFor
\end{algorithmic}
\end{algorithm}

\subsection{Algorithm}
In this work, we propose a novel variant of EXTRA, namely TT-EXTRA, to solve smooth non-convex distributed optimization problems. 
TT-EXTRA incorporates two step-sizes, allowing the construction of a Lyapunov function for the convergence analysis.

The agent view of Two-timescale EXTRA (TT-EXTRA) is presented in Algorithm \ref{alg:TTEXTRA}. Each agent $i$ holds two local variables $\mathbf{x}_{i}^r$ and $\mathbf{y}_{i}^r$. At each iteration, each agent performs local updates on $\mathbf{x}_i^r$ to obtain  $\mathbf{x}_i^{r+1}$. They, then, exchange local estimates $\mathbf{x}_{i}^{r+1}$ with their neighbors. Subsequently, the agents update $\mathbf{y}_i^{r}$. 
Note that in the updating step of $\mathbf{x}^{r+1}_i$ (line 7 in algorithm~\ref{alg:TTEXTRA}), the first three terms  correspond to distributed gradient descent component while the term involving $\mathbf{y}_i^r$ serves as a correction mechanism. This correction addresses the limitation of standard distributed gradient descent (DGD) with constant step-size, which typically converges only to a neighborhood of an optimal solution. In TT-EXTRA, $\beta$ and $\rho$ are two step-size parameters, $W$ and $\Tilde{W}$ represent two mixing matrices. Later, we will introduce conditions on these two step-sizes and two mixing matrices to ensure convergence guarantees of TT-EXTRA.

\begin{remark}
    If the two step-sizes $\beta$ and $\rho$ in Algorithm~\ref{alg:TTEXTRA} satisfy $\rho / \beta = 1$ then TT-EXTRA is equivalent to EXTRA introduced in~\cite{shi2015extra}.
\end{remark}

\begin{remark}
TT-EXTRA can be interpreted as performing Gradient Descent-Ascent (GDA) on a modified augmented Lagrangian function; see more details in Section \ref{sec:convergence}.
\end{remark}

\begin{remark}
Similar to the original EXTRA algorithm~\cite{shi2015extra}, TT-EXTRA requires the exchange of only one variable per iteration, while the variable  $\mathbf{y}_{i}$ is updated using local information.
\end{remark}
\subsection{Assumptions}
In this subsection, we state the assumptions on the objective function $f$ in~\eqref{P2} and mixing matrices $W$ and $\Tilde{W}$ that are required by the convergence analysis for TT-EXTRA.

\begin{assumption} \label{AssF}
    The objective function $f$ in Problem~\eqref{P2} is differentiable and $l$-smooth, i.e. there exists a constant $l > 0$ such that $\|\nabla f(\mathbf{x}) - \nabla f(\mathbf{y})\| < l\|\mathbf{x}-\mathbf{y}\|, \forall \mathbf{x},\mathbf{y}$.
\end{assumption}

The $l$-smooth assumption is common in the optimization literature for convergence analysis of gradient based algorithms; see \cite{shi2015extra,hong2016decomposing,nedic2017achieving}. Note that, Assumption~\ref{AssF} follows directly from the $l_{i}$-smoothness of each $f_{i}$ by defining $l = \max\{l_{1},\dots,l_{n}\}$.

\begin{assumption} \label{assf}
    The objective function $f$ in~\eqref{P2} is non-negative, i.e. $f(\mathbf{x}) \geq 0, \forall \mathbf{x}$.
\end{assumption}

Assumption \ref{assf} is equivalent to $f$ being lower bounded. This is because if there exists constant $a$ such that $f(\mathbf{x}) \geq a, \forall \mathbf{x}$, we can rewrite the optimization problem by using non-negative function $f(\mathbf{x})-a$ without changing the optimizer. Several frequently used objective functions in machine learning that satisfy Assumption \ref{AssF} and \ref{assf} are listed in \cite{hong2016decomposing}.

\begin{assumption} \label{AssM}
    Consider a connected communication graph $\mathcal{G}=(\mathcal{V},\mathcal{E})$ with set of vertices $\mathcal{V}=\{1,\dots,n\}$ and set of edges $\mathcal{E}\subseteq\mathcal{V}\times\mathcal{V}$. The two mixing matrices $W$ and $\Tilde{W}$ satisfy:
    \begin{itemize}
        \item (Decentralized property) If $(i,j) \notin \mathcal{E} $ and $i \neq j$, then $W_{ij} = \Tilde{W}_{ij} = 0$;
        \item (Symmetry) $W = W^T$, $\Tilde{W} = \Tilde{W}^T$;
        \item (Null space property) $\nullop\{W \!-\! \Tilde{W} \} \!=\! \spanop\{\bm{1}\}$ and $\nullop \{I_{n} \!-\! \Tilde{W} \} \!\supseteq \!\spanop\{\bm{1}\}$;
        \item (Spectral property) $\frac{I_{n} + W}{2} \succeq \Tilde{W} \succeq \frac{I_{n}+(\frac{1}{\rho}+1)W}{\frac{1}{\rho} + 2}$,
    \end{itemize}
    where $\rho$ is a step-size in Algorithm~\ref{alg:TTEXTRA}.
\end{assumption}

Assumption~\ref{AssM} directly implies that $\nullop(I_{n}-W) = \spanop\{\bm 1\}$ (see Prop 2.2 in \cite{shi2015extra}). Comparing Assumption~\ref{AssM} to the assumptions for EXTRA in convex settings~\cite{shi2015extra}, we tighten the range of $\Tilde{W}$ but remove the condition that $\Tilde{W}$ is positive definite. Any matrix $\Tilde{W}$ within this range that satisfies the null space properties is admissible and in our numerical result we demonstrate that it is possible to improve the convergence performance by appropriately selecting $W$ and $\Tilde{W}$. Note that the new lower bound for $\Tilde{W}$, i.e. ${(I_{n}+(\frac{1}{\rho}+1)W)}/({\frac{1}{\rho} + 2})$, is monotonically increasing in $\rho$. As $\rho \rightarrow 0$, the lower bound tends to $W$ while, as $\rho \rightarrow \infty$, it tends to the upper bound ${(I_{n} + W)}/{2}$. 

\subsection{Convergence analysis for TT-EXTRA}

The following theorem presents conditions on $\rho$ and $\beta$, the step-sizes in Algorithm~\ref{alg:TTEXTRA}, to ensure the convergence of TT-EXTRA. 

\begin{thm} \label{Thm}
     Consider the problem in (\ref{P2}). Let Assumption \ref{AssF}, \ref{assf} and \ref{AssM} hold. Let $A^T A = (\Tilde{W}-W)\otimes I_{p}$. If
     \begin{subequations}
         \begin{align}
         \rho &\!>\! 1 + \lambda_{\max} (\Tilde{W}-W),\label{conr1} \\ \rho &\!>\! \frac{8l + \sqrt{64l^2 + 16l(1-\lambda_{2}(W))}}{2(1-\lambda_{2}(W))},\label{conr2} \\ 
         \beta &\!>\! (\rho+1)\lambda_{\max}(\Tilde{W}-W)+1, \label{conb1}\\
         \beta &\!>\! \frac{\frac{L}{2}\! + \!\frac{1+2\rho}{\rho^2 (1-\lambda_{2}(W))} \!\!\left(\!4l^2\! +\! 4\rho^2 \|I_{n}\!-\!\Tilde{W}\!-\!(\Tilde{W}\!-\!W)\|^2 \!\right)}{1-1/a},\label{conb2}
    \end{align} 
     \end{subequations}
    where $L = l + \rho \|I_{n}-\Tilde{W}\|$ and $1 < a \le \frac{\rho^2 (1-\lambda_{2}(W))}{4l(1+2\rho)}$. Then, the following statements hold. 
    \begin{itemize}
        \item (Convergence) The iterate $\mathbf{x}^r= (\mathbf{x}_{1}^r,\mathbf{x}_{2}^r,\dots, \mathbf{x}_{n}^r)^T$ generated by TT-EXTRA converges to the set of first-order consensual stationary points $\mathcal{X}^*$.
        \item (Sub-linear rate) Let $T_{\epsilon} = \min\{r: \|\nabla f(\mathbf{x}^{r})+A^T \lambda^{r}+((I_{n}-\Tilde{W})\otimes I_{p})\mathbf{x}^{r}\|^2 + \|A\mathbf{x}^{r}\|^2 \le \epsilon \}$. We have $T_{\epsilon} = \mathcal{O}(\epsilon^{-1})$.
    \end{itemize}
\end{thm}

Note that the upper bound on $a$ is always greater than~$1$, i.e., $\rho^2 (1-\lambda_{2}(W))/4l(1+2\rho) > 1$, if $\rho$ satisfies \eqref{conr2}. This result establishes that TT-EXTRA will asymptotically achieve both consensus and stationarity across all agents for non-convex distributed optimization problems. Notably, the specified threshold of $\rho$ in Theorem \ref{Thm} only depends on $W$ as $\lambda_{\max}(\Tilde{W}-W) \le 1$ and the lower bound of $\beta$ relies on $W$, $\Tilde{W}$ and $\rho$. Moreover, the bounds on $\Tilde{W}$ are functions of both $\rho$ and $W$, as outlined in Assumption~\ref{AssM}. This hierarchical dependency allows for a sequential parameter selection strategy, which will be detailed in Section~\ref{sec:para selec}. Furthermore, under Assumptions~\ref{AssM}, $\|\nabla f(\mathbf{x}^{r})+A^T \lambda^{r}+(I_{np}-\Tilde{W}\otimes I_{p})\mathbf{x}^{r}\|^2 + \|A\mathbf{x}^{r}\|^2 = 0$ is a sufficient condition for Definition~\ref{def}. Thus we conclude that TT-EXTRA converges to a first-order consensual stationary point sub-linearly.

\section{Proof Sketches for Theorem~\ref{Thm}} 
\label{sec:convergence}
In this section, we present the sketch of the proof for Theorem \ref{Thm} due to space constraints. The complete proof can be found in~\cite{peng2024two}. 

First, we show that TT-EXTRA can be reformulated as:
\begin{subequations}\label{MEXTRA}
    \begin{align}
    \mathbf{x}^{r+1} \!&=\! \mathbf{x}^{r} \!-\! \frac{1}{\beta} (\nabla f(\mathbf{x}^{r}) \!+\! A^T \lambda^{r} \!+\! \rho (I_{np}\!-\!\Tilde{W}\otimes I_{p}) \mathbf{x}^{r}),  \\
    \lambda^{r+1} \!&=\! \lambda^{r} + \rho A \mathbf{x}^{r+1}, 
    \end{align}
\end{subequations}
where matrix $A$ satisfies $A^T A = (\Tilde{W}-W$)$\otimes I_{p}$ and $\lambda^{r}$ satisfies $\mathbf{y}^{r}=A^T \lambda^{r}$. Thus, TT-EXTRA can be viewed as a gradient decent step in the primal space and a gradient ascent step in the dual space. Especially, TT-EXTRA can be interpreted as the Gradient Descent Ascent method (GDA) applied to the following augmented Lagrangian-like function:
\begin{align}
    L_{\rho}(\mathbf{x},\lambda) = f(\mathbf{x}) + \langle \lambda, A\mathbf{x} \rangle + \frac{\rho}{2} \|\mathbf{x}\|_{(I_{n}-\Tilde{W})\otimes I_{p}}^2 ,\label{primal-obj}
\end{align}
where $1/\beta$ is the primal step-size and  $\rho$ is the dual step-size. 
This type of gradient decent and ascent algorithms have been studied in \cite{lu2021linearized,hong2018gradient,10035454}. However, the proofs in these papers cannot be directly applied to our problem and algorithm, primarily due to the presence of two distinct mixing matrices $W$ and $\Tilde{W}$. This structural difference leads to the $\mathbf{x}$-update in \eqref{MEXTRA} that does not correspond to a standard gradient descent step on the conventional augmented Lagrangian function. Note that when $\Tilde{W}=(I_{n}+W)/2$, $L_{\rho}(\mathbf{x},\lambda)$ defined in \eqref{primal-obj} reduces to the standard augmented Lagrangian function.

\subsection{Key lemmas}
In the following, we establish the convergence analysis of TT-EXTRA and provide a sketch of the proof of Theorem \ref{Thm} based on (\ref{MEXTRA}). Let $C_{1} = \beta I_{np} - \rho A^T A$ ($C_{1} \succeq 0$ if $\beta$ satisfies \eqref{conb1}) and $ C_{2} = \rho (I_{n}-\Tilde{W} -(\Tilde{W}-W) )\otimes I_{p}$ ($C_{2} \succeq 0$ if $W,\Tilde{W}$ satisfy Assumption~\ref{AssM}). We construct the following Lyapunov function, which will play a key role in our analysis:
\begin{align}
     P&_{c}(\mathbf{x}^{r+1},\mathbf{x}^{r},\lambda^{r+1})  = L_{\rho}(\mathbf{x}^{r+1},\lambda^{r+1})+ \nonumber \\  &\!+\! \frac{c\rho}{2} \|A\mathbf{x}^{r+1}\|^2 + \frac{c}{2} \|\mathbf{x}^{r+1}-\mathbf{x}^{r}\|^2_{C_{1}+C_{2}} \nonumber \\ &\!+\! \frac{4 \rho^2(1\!+\!2\rho)}{\rho^2 (1\!-\!\lambda_{2}(W))}  \|(I_{n}\!-\!2\Tilde{W}\!+W)\otimes I_{p}\|^2 \|\mathbf{x}^{r+1}\!-\!\mathbf{x}^{r}\|^2 \nonumber\\ 
     & \!+\! (\frac{cl}{2}+\frac{4l^2(1\!+\!2\rho)}{\rho^2 (1\!-\!\lambda_{2}(W))})  \|\mathbf{x}^{r+1}-\mathbf{x}^{r}\|^2. \label{lyafun}
\end{align}
where $c$ denotes a positive constant depending on $\beta$ (see the proof of Lemma~\ref{Le4} in \cite{peng2024two}).

We then present the following five lemmas that will be used in the proof of Theorem~\ref{Thm}. These lemmas particularly show that the Lyapunov function defined in \eqref{lyafun} is monotonically decreasing and lower bounded. The proofs of these lemmas can be found in~\cite{peng2024two}.

\begin{lemma} \label{Le1}
Let $\mathbf{x}^{r},\lambda^{r}$ be the sequence generated by \eqref{MEXTRA}. If Assumption \ref{AssF} holds, we have
\begin{align} \label{lemmaL}
     L_{\rho}&(\mathbf{x}^{r+1},\lambda^{r+1}) - L_{\rho}(\mathbf{x}^{r},\lambda^{r})  \nonumber \\ & \le \frac{1}{\rho} \|\lambda^{r+1}-\lambda^{r}\|^2 - (\beta - \frac{L}{2}) \|\mathbf{x}^{r+1}-\mathbf{x}^{r}\|^2 
\end{align}
where $L_{\rho}(.,.)$ is defined in \eqref{primal-obj} and $L = l + \rho \|I_{n}-\Tilde{W}\|$.
\end{lemma}

\begin{lemma} \label{Le2}
   Let $\lambda_{2}(W) < 1$ be the second largest eigenvalue of $W$, $\mathbf{x}^{r},\lambda^{r}$ be the sequence generated by \eqref{MEXTRA} and $\mathbf{w}^{r+1} = (\mathbf{x}^{r+1}-\mathbf{x}^{r}) - (\mathbf{x}^{r}- \mathbf{x}^{r-1})$. If Assumptions \ref{AssF}, \ref{assf} and \ref{AssM} hold, we have 
\begin{align}
    & \ \frac{1}{\rho}\|\lambda^{r+1}-\lambda^{r}\|^2 \nonumber \\ \le &\frac{1+2\rho}{\rho^2 (1-\lambda_{2}(W))} \left(4l^2 \|\mathbf{x}^{r}-\mathbf{x}^{r-1}\|^2 + 2\|\mathbf{w}^{r+1}\|_{C_{1}^T C_{1}}^2\right) \nonumber \\ &\!+\! \frac{4\rho^2 (1\!+\!2\rho)}{\rho^2 (1\!-\!\lambda_{2}(W))}\|(I_{n}\!\!-\!2\Tilde{W}\!+\!W)\otimes I_{p}\|^2 \|\mathbf{x}^{r}\!-\!\mathbf{x}^{r-1}\|^2.\label{lemma2}
\end{align} 
\end{lemma}

\begin{lemma} \label{Le3}
   Let $C_{3}=\beta I_{np} - \rho (I_{n}-\Tilde{W})\otimes I_{p}$ and  $\mathbf{x}^{r},\lambda^{r}$ be the sequence generated by \eqref{MEXTRA}. If Assumptions \ref{AssF},\ref{assf} and \ref{AssM} hold, we have
    \begin{align}
        &\frac{\rho}{2} \|A\mathbf{x}^{r+1}\|^2 + \frac{1}{2} \|\mathbf{x}^{r+1}-\mathbf{x}^{r}\|^2_{C_{1}} + \frac{1}{2}\|\mathbf{x}^{r+1}-\mathbf{x}^{r}\|^2_{C_{2}} \nonumber \\
        \le \ &\frac{\rho}{2} \|A\mathbf{x}^{r}\|^2 + \frac{1}{2} \|\mathbf{x}^{r}-\mathbf{x}^{r-1}\|^2_{C_{1}} + \frac{1}{2}\|\mathbf{x}^{r}-\mathbf{x}^{r-1}\|^2_{C_{2}} + \nonumber \\
        & \frac{l}{2} \|\mathbf{x}^{r}\!-\!\mathbf{x}^{r-1}\|^2 \!+\! \frac{l}{2} \|\mathbf{x}^{r+1}\!-\!\mathbf{x}^{r}\|^2 \!-\! \frac{1}{2}\|\mathbf{w}^{r+1}\|^2_{C_{3}}. \label{lemma3}
    \end{align}
\end{lemma}

\begin{lemma} \label{Le4}
    Let $\mathbf{x}^{r},\lambda^{r}$ be the sequence generated by \eqref{MEXTRA}. If Assumptions \ref{AssF}, \ref{assf} and \ref{AssM} hold, $\rho$ and $\beta$ satisfy (\ref{conr1}), (\ref{conr2}), (\ref{conb1}) and (\ref{conb2}), then the Lyapunov function $P_{c}(\mathbf{x}^{r+1},\mathbf{x}^{r},\lambda^{r+1})$ defined in (\ref{lyafun}) decreases monotonically. \label{lemma4}
\end{lemma}

\begin{lemma} \label{Le5}
    Let $\mathbf{x}^{r},\lambda^{r}$ be the sequence generated by \eqref{MEXTRA}. If  Assumptions \ref{AssF},\ref{assf} and \ref{AssM} hold, then there exists a constant \underline{$p$} $< \infty$ such that $P_{c} (\mathbf{x}^{r+1},\mathbf{x}^{r},\lambda^{r+1}) \geq $ \underline{$p$} $, \forall r \geq 0$. \label{lemma5}
\end{lemma}
\subsection{Proof Sketch for Theorem~\ref{Thm}}
Combining Lemmas~\ref{Le1}, \ref{Le2}, and \ref{Le3} results in
\begin{align}
    P_{c} (\mathbf{x}^{r+1}&,\mathbf{x}^{r},\lambda^{r+1})- P_{c} (\mathbf{x}^{r},\mathbf{x}^{r-1},\lambda^{r}) \nonumber\\ \le &- \|\mathbf{w}^{r+1}\|_{B}^2 -b\|\mathbf{x}^{r+1}-\mathbf{x}^{r}\|^2 \label{descent}
\end{align}
where $B=\frac{c}{2}C_{3} - \frac{2(1+2\rho)}{\rho^2 (1-\lambda_{2}(W))}C_{1}^T C_{1}$ and  $b=\beta - cl - \frac{L}{2} - \frac{1+2\rho}{\rho^2 (1-\lambda_{2}(W))} (4l^2 + 4\rho^2 \|(I_{n}-2\Tilde{W}+W)\otimes I_{p}\|^2 )$. Lemma~\ref{Le4} shows that the matrix $B$ is positive definite and the constant $b$ is positive. Furthermore, Lemma~\ref{Le5} implies that $P_{c} (\mathbf{x}^{r+1},\mathbf{x}^{r},\lambda^{r+1})$ is lower bounded. Thus, we establish the existence of the limit of \( P_{c}(\mathbf{x}^{r+1}, \mathbf{x}^{r}, \lambda^{r+1}) \), which, together with the descent property in~\eqref{descent}, implies that $\lim_{r \rightarrow \infty} \| \mathbf{x}^{r+1} - \mathbf{x}^{r} \| = 0$. This result, when combined with Lemma~\ref{Le2}, ensures asymptotic consensus, and when further combined with the $\mathbf{x}$-update in~\eqref{MEXTRA}, guarantees the asymptotic stationarity. The sub-linear convergence rate, i.e. $T_{\epsilon} = \mathcal{O}(1/\epsilon)$, can subsequently be derived using the descent of the Lyapunov function discussed in~\eqref{descent}. The detailed proof of Theorem \ref{Thm} can be found in~\cite{peng2024two}.

\section{Parameter Selection} \label{sec:para selec}
In this section, we introduce a sequential method for selecting the two mixing matrices $W$ and $\Tilde{W}$ and the two step-sizes $\beta$ and $\rho$ to ensure they meet the conditions established in Assumption~\ref{AssM} and Theorem~\ref{Thm}. Note that, to guarantee the convergence of TT-EXTRA, the mixing matrix \( W \) is only required to satisfy the null space conditions specified in Assumption~\ref{AssM}. The lower bounds of \( \rho \) in~\eqref{conr1} and~\eqref{conr2} depend solely on \( W \). Subsequently, the spectral properties of \( \Tilde{W} \) in Assumption~\ref{AssM} are determined by both \( W \) and \( \rho \). Finally, the lower bound of \( \beta \) is influenced by \( W \), \( \Tilde{W} \), and \( \rho \). This hierarchical dependency enables for sequential selection of the parameters. The detailed steps are presented in Algorithm~\ref{paras}.

\begin{algorithm}[H]
\caption{Parameter Selection}\label{paras}
\begin{algorithmic}[1]
    \State Input: a connected network  $\mathcal{G}$
    \State Set $W$ based on $\mathcal{G}$ such that $\nullop(I_{n}-W) = \spanop\{\bm 1 \}$
    \State Set $\rho >$ \underline{$\rho$}, where lower bound  \underline{$\rho$} satisfies~\eqref{conr1} and \eqref{conr2}
    \State Set $\Tilde{W} = \frac{I_{n}+(\frac{1}{\rho}+1)W}{\frac{1}{\rho} + 2}$ with $\rho >$ \underline{$\rho$}
    \State Choose $\beta$ satisfying \eqref{conb1} and \eqref{conb2}
    \State Output: $(W,\rho,\Tilde{W},\beta)$
\end{algorithmic}
\end{algorithm}
\begin{remark}
    The implementation of Algorithm~\ref{paras} requires knowledge of the global network, as the lower bounds of $\rho$ and $\beta$ are determined by the eigenvalues of the mixing matrices; see \eqref{conr2} and \eqref{conb2}.
\end{remark}

Note that in line 1, some feasible and frequently used $W$ such as Laplacian-based constant edge weight matrix \cite{sayed2014diffusion} and Metropolis constant edge weight matrix \cite{boyd2004fastest} are listed in \cite{shi2015extra}. In line 4, the choice of $\Tilde{W}$ is not unique and can be any matrix of the form $\Tilde{W} = \frac{I_{n}+(\frac{1}{\rho}+1)W}{\frac{1}{\rho} + 2}$ with $\rho >$ \underline{$\rho$}. The following two lemmas which indicate that $\Tilde{W}$ satisfies the null space and spectral properties stated in Assumption~\ref{AssM}.
\begin{lemma} \label{pro}
 Let Assumption~\ref{AssM} holds. If $\rho_{1} \geq \rho_{2}>0$, then $\frac{I_{n}+(\frac{1}{\rho_{1}}+1)W}{\frac{1}{\rho_{1}}+2} \succeq \frac{I_{n} + (\frac{1}{\rho_{2}}+1)W}{\frac{1}{\rho_{2}}+2}$.
\end{lemma}

\begin{lemma} \label{prop3}
    Let Assumption \ref{AssM} holds, and $\Tilde{W} = aI_{n} + bW$ for $a,b \in \mathbb{R}$. If $a+b = 1$ and $\nullop(I_{n}-W) = \spanop\{\bm 1 \}$ then we have $\nullop(\Tilde{W}-W) = \spanop \{\bm 1 \}$ and $\nullop(I_{n}-\Tilde{W}) \supseteq \spanop \{\bm 1 \}$. 
\end{lemma}

The complete proof of these two lemmas can be found in~\cite{peng2024two} and are omitted here due to space constraints. The following proposition states that the output of Algorithm \ref{paras} always satisfies the conditions outlined in the previous section to ensure the convergence of TT-EXTRA.

\begin{prop} \label{para selec}
    The parameters $(W,\rho,\Tilde{W},\beta)$ generated by Algorithm~\ref{paras} satisfy the conditions in Assumption \ref{AssM} and Theorem~\ref{Thm}.
\end{prop} 

\begin{proof}
    In line 3 of Algorithm \ref{paras}, we may set \underline{$\rho$} $= max\{1 + \lambda_{\max} (\Tilde{W}-W), \frac{8l + \sqrt{64l^2 + 16l(1-\lambda_{2}(W))}}{2(1-\lambda_{2}(W))} \}$. By choosing \(\rho > \underline{\rho}\), the conditions in \eqref{conr1} and \eqref{conr2} are satisfied. In line 4 of Algorithm \ref{paras}, we set \( \Tilde{W} = \big(I_{n} + (\frac{1}{\rho} + 1)W\big) / \big(\frac{1}{\rho} + 2\big) \) with $\rho >$ \underline{$\rho$}. Combined with the monotonicity property established in Lemma~\ref{pro}, this implies that $\Tilde{W}$ satisfies the spectral conditions specified in Assumption~\ref{AssM}. Additionally, in  line 1 of Algorithm~\ref{paras}, we select $W$ to satisfy $\nullop(I_{n}-W) = \spanop\{\bm 1 \}$. Together with the fact that $\Tilde{W}$ satisfies the conditions stated in Lemma~\ref{prop3} for any \( \rho > 0 \), it implies that $W$ and $\Tilde{W}$ satisfy the null space properties specified in Assumption~\ref{AssM}. Therefore, we conclude that the parameters $(W,\rho,\Tilde{W},\beta)$ generated by Algorithm~\ref{paras} satisfy the conditions required in Assumption \ref{AssM} and Theorem~\ref{Thm}.
\end{proof}

\begin{figure}[h]
  \includegraphics[width=1\linewidth]{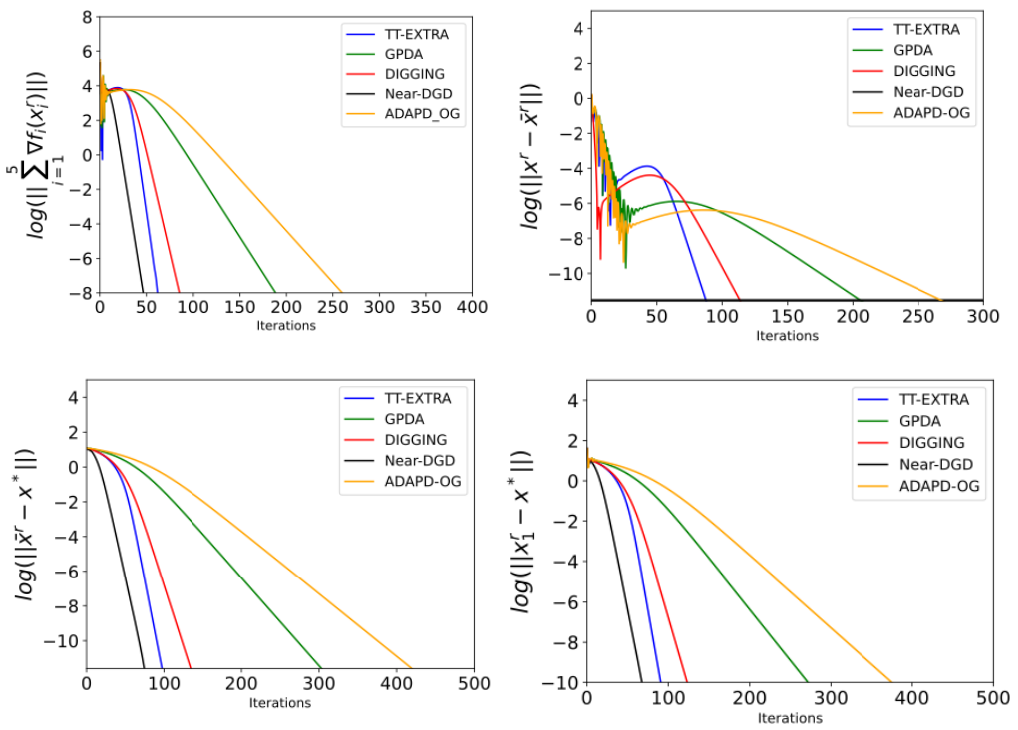}

\caption{Comparison of TT-EXTRA with four algorithms: GPDA~\cite{hong2018gradient}, DIGGING~\cite{nedic2017achieving}, ADAPD-OG~\cite{10035454} and Near-DGD~\cite{iakovidou2021convergence}.}
\label{fig:exp}
\end{figure}

\section{Numeric results} \label{sec:numeric}
We demonstrate the performance of TT-EXTRA on the following decentralized optimization problem:
\begin{align*}
    \min_{x\in \mathcal{R}} f(\mathbf{x})=\sum_{i=1}^{5} f_{i}(\mathbf{x}) 
\end{align*}
where for $i \in \{1,\dots, 5 \}$
\begin{align}
    f_{i}(\mathbf{x}) &\!=\! \begin{cases} 
       a_{1}\mathbf{x}^4 \!+\!a_{2}\mathbf{x}^3 \!+\! a_{3}\mathbf{x}^2 \!+\! a_{4}\mathbf{x} & |\mathbf{x}| \!\leq\! 10 \\
      b_{1}\mathbf{x} \!-\! b_{2} & \mathbf{x} \!<\! -10 \\
      c_{1}\mathbf{x} \!-\! c_{2} & \mathbf{x} \!>\! 10
   \end{cases}
   \label{exp-obj} 
\end{align}
\begin{table}[h] 
\centering
\begin{tabular}{|c|c|c|c|c|c|c|c|c|}
\hline
&$a_{1}$ & $a_{2}$ & $a_{3}$ & $a_{4}$ & $b_{1}$ & $b_{2}$ & $c_{1}$ &$c_{2}$ \\
\hline
$f_{1}$&1 & -4 & 0 & 0 & -5200 & -38000 & 2800 & -22000 \\
\hline
$f_{2}$&0.5 & 0 & -3 & 0 & -1940 & -14700 & 1940 & -14700 \\
\hline
$f_{3}$&-0.5 & 2 & -4 & 0 & 2680 & 19400 & -1480 & 11400 \\
\hline
$f_{4}$&0.5 & -1 & 0 & 3 & -2297 & -17000 & 1703 & -13000 \\
\hline
$f_{5}$&-1 & 0 & 5 & -7 &3893 & 29470 & -3907 & 29570 \\
\hline
\end{tabular}
\caption{Polynomial coefficients}
\label{tab:exp}
\end{table}

The coefficients of each local objective function $f_i$ in \eqref{exp-obj} are provided in Table~\ref{tab:exp}. In this setup, each $f_i$ is continuously differentiable, and the overall objective function $f$ is $l$-smooth with $l = 616$, and admits both a saddle point and a global minimizer denoted by $\mathbf{x}^*$. Similar objective functions have been considered in the literature including~\cite{tatarenko2017non,zeng2018nonconvex}.

The connected network with 5 nodes is randomly generated and the mixing matrix $W$ was selected to be the Laplacian-based constant edge weight matrix\cite{sayed2014diffusion}. We compare our TT-EXTRA with other four distributed algorithms including GPDA \cite{hong2018gradient}, DIGGING \cite{nedic2017achieving}, Near-DGD (with 50 consensus rounds per gradient evaluation)~\cite{iakovidou2021convergence} and ADAPD-OG~\cite{10035454}. All algorithms were initialized from the same randomly selected starting point. The stepsizes for each algorithm are manually tuned to achieve the largest values that ensure convergence, in order to evaluate their best achievable performance for fairness.

Among the five algorithms, TT-EXTRA, GPDA, and ADAPD-OG require the exchange of only one local variable per iteration, whereas DIGGING involves the exchange of two. Near-DGD performs 50 exchanges of a single local variable. In this experiment, performing 50 consensus rounds per gradient evaluation is the minimum required to guarantee the robust convergence of Near-DGD.


Figure~\ref{fig:exp} presents the performance of the five algorithms regarding  the decrease of gradient, consensus error, and objective error. Near-DGD demonstrates the fastest convergence among all algorithms. However, this rapid convergence comes at the higher communication cost caused by the repeated execution of consensus steps per iteration. Due to the flexibility in the design of the two stepsizes and the matrices $W$ and $\Tilde{W}$, this example demonstrates that TT-EXTRA converges faster than other algorithms except for Near-DGD. Overall, we observe that TT-EXTRA achieves a favorable balance between communication efficiency and convergence speed.


\section{CONCLUSIONS}\label{sec:conc}
In conclusion, this paper presents a convergence analysis for a variant of EXTRA, referred to as two timescale EXTRA (TT-EXTRA), for smooth non-convex distributed optimization problems. We establish the sub-linear convergence of TT-EXTRA to first-order consensual stationary points. We also introduce a sequential approach for parameter selection. Future work can focus on using noise so that TT-EXTRA can escape from the saddle point within polynomial number of iterations using the developed Lyapunov function.

\bibliographystyle{plain}
\bibliography{sample}

\begin{thebibliography}{10}

\bibitem{alghunaim2020linear}
Sulaiman~A Alghunaim and Ali~H Sayed.
\newblock Linear convergence of primal--dual gradient methods and their performance in distributed optimization.
\newblock {\em Automatica}, 117:109003, 2020.

\bibitem{alghunaim2022unified}
Sulaiman~A Alghunaim and Kun Yuan.
\newblock A unified and refined convergence analysis for non-convex decentralized learning.
\newblock {\em IEEE Transactions on Signal Processing}, 70:3264--3279, 2022.

\bibitem{boyd2004fastest}
Stephen Boyd, Persi Diaconis, and Lin Xiao.
\newblock Fastest mixing markov chain on a graph.
\newblock {\em SIAM review}, 46(4):667--689, 2004.

\bibitem{chang2014multi}
Tsung-Hui Chang, Mingyi Hong, and Xiangfeng Wang.
\newblock Multi-agent distributed optimization via inexact consensus {ADMM}.
\newblock {\em IEEE Transactions on Signal Processing}, 63(2):482--497, 2014.

\bibitem{daneshmand2020second}
Amir Daneshmand, Gesualdo Scutari, and Vyacheslav Kungurtsev.
\newblock Second-order guarantees of distributed gradient algorithms.
\newblock {\em SIAM Journal on Optimization}, 30(4):3029--3068, 2020.

\bibitem{di2015distributed}
Paolo Di~Lorenzo and Gesualdo Scutari.
\newblock Distributed nonconvex optimization over networks.
\newblock In {\em 2015 IEEE 6th International Workshop on Computational Advances in Multi-Sensor Adaptive Processing (CAMSAP)}, pages 229--232. IEEE, 2015.

\bibitem{di2016distributed}
Paolo Di~Lorenzo and Gesualdo Scutari.
\newblock Distributed nonconvex optimization over time-varying networks.
\newblock In {\em 2016 IEEE International Conference on Acoustics, Speech and Signal Processing (ICASSP)}, pages 4124--4128. IEEE, 2016.

\bibitem{facchinei2015parallel}
Francisco Facchinei, Gesualdo Scutari, and Simone Sagratella.
\newblock Parallel selective algorithms for nonconvex big data optimization.
\newblock {\em IEEE Transactions on Signal Processing}, 63(7):1874--1889, 2015.

\bibitem{gan2012optimal}
Lingwen Gan, Ufuk Topcu, and Steven~H Low.
\newblock Optimal decentralized protocol for electric vehicle charging.
\newblock {\em IEEE Transactions on Power Systems}, 28(2):940--951, 2012.

\bibitem{gao2019reinforcement}
Weinan Gao, Jingqin Gao, Kaan Ozbay, and Zhong-Ping Jiang.
\newblock Reinforcement-learning-based cooperative adaptive cruise control of buses in the lincoln tunnel corridor with time-varying topology.
\newblock {\em IEEE Transactions on Intelligent Transportation Systems}, 20(10):3796--3805, 2019.

\bibitem{hong2016decomposing}
Mingyi Hong.
\newblock Decomposing linearly constrained nonconvex problems by a proximal primal dual approach: Algorithms, convergence, and applications.
\newblock {\em arXiv preprint arXiv:1604.00543}, 2016.

\bibitem{hong2017prox}
Mingyi Hong, Davood Hajinezhad, and Ming-Min Zhao.
\newblock Prox-pda: The proximal primal-dual algorithm for fast distributed nonconvex optimization and learning over networks.
\newblock In {\em International Conference on Machine Learning}, pages 1529--1538. PMLR, 2017.

\bibitem{hong2018gradient}
Mingyi Hong, Meisam Razaviyayn, and Jason Lee.
\newblock Gradient primal-dual algorithm converges to second-order stationary solution for nonconvex distributed optimization over networks.
\newblock In {\em International Conference on Machine Learning}, pages 2009--2018. PMLR, 2018.

\bibitem{iakovidou2021convergence}
Charikleia Iakovidou and Ermin Wei.
\newblock On the convergence of near-dgd for nonconvex optimization with second order guarantees.
\newblock In {\em 2021 60th IEEE Conference on Decision and Control (CDC)}, pages 259--264. IEEE, 2021.

\bibitem{ling2015dlm}
Qing Ling, Wei Shi, Gang Wu, and Alejandro Ribeiro.
\newblock {DLM}: Decentralized linearized alternating direction method of multipliers.
\newblock {\em IEEE Transactions on Signal Processing}, 63(15):4051--4064, 2015.

\bibitem{lu2021linearized}
Songtao Lu, Jason~D Lee, Meisam Razaviyayn, and Mingyi Hong.
\newblock Linearized {ADMM} converges to second-order stationary points for non-convex problems.
\newblock {\em IEEE Transactions on Signal Processing}, 69:4859--4874, 2021.

\bibitem{10035454}
Gabriel Mancino-Ball, Yangyang Xu, and Jie Chen.
\newblock A decentralized primal-dual framework for non-convex smooth consensus optimization.
\newblock {\em IEEE Transactions on Signal Processing}, 71:525--538, 2023.

\bibitem{mokhtari2015decentralized}
Aryan Mokhtari and Alejandro Ribeiro.
\newblock Decentralized double stochastic averaging gradient.
\newblock In {\em ACSSC}, pages 406--410, 2015.

\bibitem{mokhtari2016dqm}
Aryan Mokhtari, Wei Shi, Qing Ling, and Alejandro Ribeiro.
\newblock Dqm: Decentralized quadratically approximated alternating direction method of multipliers.
\newblock {\em IEEE Transactions on Signal Processing}, 64(19):5158--5173, 2016.

\bibitem{nedic2017achieving}
Angelia Nedic, Alex Olshevsky, and Wei Shi.
\newblock Achieving geometric convergence for distributed optimization over time-varying graphs.
\newblock {\em SIAM Journal on Optimization}, 27(4):2597--2633, 2017.

\bibitem{nedic2009distributed}
Angelia Nedic and Asuman Ozdaglar.
\newblock Distributed subgradient methods for multi-agent optimization.
\newblock {\em IEEE Transactions on Automatic Control}, 54(1):48--61, 2009.

\bibitem{peng2024two}
Zeyu Peng, Farhad Farokhi, and Ye~Pu.
\newblock Two timescale extra for smooth non-convex distributed optimization problems.
\newblock {\em arXiv preprint arXiv:2411.19483}, 2024.

\bibitem{ram2009distributed}
Sundhar~Srinivasan Ram, Venugopal~V Veeravalli, and Angelia Nedic.
\newblock Distributed non-autonomous power control through distributed convex optimization.
\newblock In {\em IEEE INFOCOM 2009}, pages 3001--3005. IEEE, 2009.

\bibitem{sayed2014diffusion}
Ali~H Sayed.
\newblock Diffusion adaptation over networks.
\newblock In {\em Academic Press Library in Signal Processing}, volume~3, pages 323--453. Elsevier, 2014.

\bibitem{schizas2007consensus}
Ioannis~D Schizas, Alejandro Ribeiro, and Georgios~B Giannakis.
\newblock Consensus in ad hoc wsns with noisy links—part i: Distributed estimation of deterministic signals.
\newblock {\em IEEE Transactions on Signal Processing}, 56(1):350--364, 2007.

\bibitem{shi2015extra}
Wei Shi, Qing Ling, Gang Wu, and Wotao Yin.
\newblock {EXTRA}: An exact first-order algorithm for decentralized consensus optimization.
\newblock {\em SIAM Journal on Optimization}, 25(2):944--966, 2015.

\bibitem{shi2014linear}
Wei Shi, Qing Ling, Kun Yuan, Gang Wu, and Wotao Yin.
\newblock On the linear convergence of the {ADMM} in decentralized consensus optimization.
\newblock {\em IEEE Transactions on Signal Processing}, 62(7):1750--1761, 2014.

\bibitem{sun2019distributed}
Haoran Sun and Mingyi Hong.
\newblock Distributed non-convex first-order optimization and information processing: Lower complexity bounds and rate optimal algorithms.
\newblock {\em IEEE Transactions on Signal processing}, 67(22):5912--5928, 2019.

\bibitem{tatarenko2017non}
Tatiana Tatarenko and Behrouz Touri.
\newblock Non-convex distributed optimization.
\newblock {\em IEEE Transactions on Automatic Control}, 62(8):3744--3757, 2017.

\bibitem{wang2019global}
Yu~Wang, Wotao Yin, and Jinshan Zeng.
\newblock Global convergence of {{ADMM}} in nonconvex nonsmooth optimization.
\newblock {\em Journal of Scientific Computing}, 78:29--63, 2019.

\bibitem{xiao2004fast}
Lin Xiao and Stephen Boyd.
\newblock Fast linear iterations for distributed averaging.
\newblock {\em Systems \& Control Letters}, 53(1):65--78, 2004.

\bibitem{yao2018distributed}
Lisha Yao, Ye~Yuan, Shreyas Sundaram, and Tao Yang.
\newblock Distributed finite-time optimization.
\newblock In {\em 2018 IEEE 14th International Conference on Control and Automation (ICCA)}, pages 147--154. IEEE, 2018.

\bibitem{yuan2016convergence}
Kun Yuan, Qing Ling, and Wotao Yin.
\newblock On the convergence of decentralized gradient descent.
\newblock {\em SIAM Journal on Optimization}, 26(3):1835--1854, 2016.

\bibitem{zeng2018nonconvex}
Jinshan Zeng and Wotao Yin.
\newblock On nonconvex decentralized gradient descent.
\newblock {\em IEEE Transactions on signal processing}, 66(11):2834--2848, 2018.

\end{thebibliography}
\appendix
\section{Supplementary Material}
\subsection{Proof of Lemma \ref{Le1}}
Since $L_{\rho}(.,\lambda)$ is $L$-smooth for all $\lambda$ with $L=l + \rho \|I_{n}-\Tilde{W}\|$, then we have:
\begin{align*}
    & \ \ \ \ \ L_{\rho}(\mathbf{x}^{r+1},\lambda^r)\\ &\le L_{\rho}(\mathbf{x}^{r},\lambda^r) + \langle \nabla_{\mathbf{x}} L_{\rho}(\mathbf{x}^{r},\lambda^r), \mathbf{x}^{r+1}-\mathbf{x}^r \rangle + \frac{L}{2} ||\mathbf{x}^{r+1}-\mathbf{x}^r||^2 \\
    &= L_{\rho}(\mathbf{x}^{r},\lambda^r) + \langle \nabla f(\mathbf{x}^r) + A^T \lambda^r +\rho (I_{np}-\Tilde{W}\otimes I_{p})\mathbf{x}^r , \mathbf{x}^{r+1}-\mathbf{x}^r \rangle + \frac{L}{2} ||\mathbf{x}^{r+1}-\mathbf{x}^r||^2 \\ 
    &= L_{\rho}(\mathbf{x}^{r},\lambda^r) + \langle \nabla f(\mathbf{x}^r) + A^T \lambda^r +\rho (I_{np}-\Tilde{W}\otimes I_{p})\mathbf{x}^r + \beta (\mathbf{x}^{r+1}-\mathbf{x}^r) - \beta (\mathbf{x}^{r+1}-\mathbf{x}^r), \mathbf{x}^{r+1}-\mathbf{x}^r \rangle + \\ & \ \ \ \frac{L}{2} ||\mathbf{x}^{r+1}-\mathbf{x}^r||^2 \\ 
    & = L_{\rho}(\mathbf{x}^{r},\lambda^r)   - \beta ||\mathbf{x}^{r+1}-\mathbf{x}^r||^2 + \frac{L}{2} ||\mathbf{x}^{r+1}-\mathbf{x}^r||^2 \\ 
    &= L_{\rho}(\mathbf{x}^{r},\lambda^r) - (\beta - \frac{L}{2}) ||\mathbf{x}^{r+1}-\mathbf{x}^r||^2
\end{align*}
where the third equality comes from \ref{MEXTRA}.
\begin{align*}
    & \ \ \ \ L_{\rho}(\mathbf{x}^{r+1},\lambda^{r+1}) - L_{\rho}(\mathbf{x}^{r+1},\lambda^{r}) \\ & = f(\mathbf{x}^{r+1}) + \langle \lambda^{r+1}, A\mathbf{x}^{r+1} \rangle + \frac{\rho}{2} ||{\mathbf{x}^{r+1}}||_{I_{np}-\Tilde{W}\otimes I_{p}}^2 - \left(f(\mathbf{x}^{r+1}) + \langle \lambda^{r}, A\mathbf{x}^{r+1} \rangle + \frac{\rho}{2} ||{\mathbf{x}^{r+1}}||_{I_{np}-\Tilde{W}\otimes I_{p}}^2 \right) \\
    &= \langle \lambda^{r+1}, A\mathbf{x}^{r+1} \rangle - \langle \lambda^{r}, A\mathbf{x}^{r+1} \rangle \\
    &= \langle \lambda^{r+1}-\lambda^r, \frac{1}{\rho} ( \lambda^{r+1}-\lambda^r) \rangle \\
    &= \frac{1}{\rho} ||\lambda^{r+1}-\lambda^r||^2
\end{align*}
Combine the above two inequalities together we have:
\begin{align*}
    & \ \ \ \ L_{\rho}(\mathbf{x}^{r+1},\lambda^{r+1}) - L_{\rho}(\mathbf{x}^{r},\lambda^{r}) \\
    &= L_{\rho}(\mathbf{x}^{r+1},\lambda^{r+1}) -L_{\rho}(\mathbf{x}^{r+1},\lambda^{r}) +L_{\rho}(\mathbf{x}^{r+1},\lambda^{r}) - L_{\rho}(\mathbf{x}^{r},\lambda^{r}) \\
    &\le - (\beta - \frac{L}{2}) ||\mathbf{x}^{r+1}-\mathbf{x}^r||^2 + \frac{1}{\rho} ||\lambda^{r+1}-\lambda^r||^2 
\end{align*} \hfill $\blacksquare$
\subsection{Proof of Lemma \ref{Le2}}
From \ref{MEXTRA} we have:
\begin{align}
    &\nabla f(\mathbf{x}^r)+A^T \lambda^r + \rho (I_{np}-\Tilde{W}\otimes I_{p}) \mathbf{x}^r + \rho ((\Tilde{W} - W)\otimes I_{p}) \mathbf{x}^r - \rho ((\Tilde{W} - W)\otimes I_{p}) \mathbf{x}^r + \beta (\mathbf{x}^{r+1}-\mathbf{x}^r) = 0 \nonumber \\
    \implies & \nabla f(\mathbf{x}^r)+A^T \lambda^r + \rho ((\Tilde{W} - W)\otimes I_{p}) \mathbf{x}^{r+1} + \rho ((\Tilde{W} - W)\otimes I_{p}) (\mathbf{x}^{r}-\mathbf{x}^{r+1}) + \nonumber \\  \ \ \ \ \ &\rho (I_{np}-\Tilde{W}\otimes I_{p}-((\Tilde{W} - W)\otimes I_{p})) \mathbf{x}^r + \beta (\mathbf{x}^{r+1}-\mathbf{x}^r) = 0 \nonumber \\
    \implies & \nabla f(\mathbf{x}^r)+A^T \lambda^{r+1} + \rho ((\Tilde{W} - W)\otimes I_{p}) (\mathbf{x}^{r}-\mathbf{x}^{r+1}) + \rho (I_{np}-\Tilde{W}\otimes I_{p}-((\Tilde{W} - W)\otimes I_{p})) \mathbf{x}^r + \nonumber \\ &\beta (\mathbf{x}^{r+1}-\mathbf{x}^r) = 0 \label{citer}
\end{align}
note that \ref{citer} is true $\forall r > 0$, therefore we also have:
\begin{align}
    \nabla f(\mathbf{x}^{r-1})+A^T \lambda^r +& \rho ((\Tilde{W} - W)\otimes I_{p}) (\mathbf{x}^{r-1}-\mathbf{x}^{r}) + \nonumber \\ &\rho \left(I_{np}-\Tilde{W}\otimes I_{p}-((\Tilde{W} - W)\otimes I_{p})\right) \mathbf{x}^{r-1} + \beta (\mathbf{x}^{r}-\mathbf{x}^{r-1}) = 0 \label{piter}
\end{align}
Taking the difference between \ref{citer} and \ref{piter} we get:
\begin{align*}
    \nabla f(\mathbf{x}^r) -& \nabla f(\mathbf{x}^{r-1}) + A^T (\lambda^{r+1}-\lambda^r) + \rho ((\Tilde{W} - W)\otimes I_{p}) \left\{ (\mathbf{x}^r-\mathbf{x}^{r+1}) - (\mathbf{x}^{r-1} - \mathbf{x}^{r}) \right\} +  \\  &\rho \left(I_{np}-\Tilde{W}\otimes I_{p}-((\Tilde{W} - W)\otimes I_{p})\right) (\mathbf{x}^r - \mathbf{x}^{r-1}) + \beta \left\{(\mathbf{x}^{r+1}-\mathbf{x}^r)-(\mathbf{x}^{r}-\mathbf{x}^{r-1})  \right\} = 0
\end{align*}
which implies:
\begin{align}
    \nabla f(\mathbf{x}^r) - \nabla f(\mathbf{x}^{r-1}) + A^T (\lambda^{r+1}-\lambda^r) + C_{1}\mathbf{w}^{r+1} + \rho \left(I_{np}-\Tilde{W}\otimes I_{p}-((\Tilde{W} - W)\otimes I_{p})\right) (\mathbf{x}^r - \mathbf{x}^{r-1}) = 0 \label{dualdiff}
\end{align}
let $A = ((\Tilde{W} - W)\otimes I_{p})^{1/2}$. Combined with $\bm 1^T (\lambda^{r+1}-\lambda^{r}) = 0, \forall r \geq 0$ and $\textup{Null}\{I-W\} = \textup{span} \{\bm 1$\} then we have:
\begin{align*}
    ||A^T (\lambda^{r+1}-\lambda^{r})||^2 &= (\lambda^{r+1}-\lambda^{r})^T A A^T (\lambda^{r+1}-\lambda^{r}) \\&= (\lambda^{r+1}-\lambda^{r})^T A^T A (\lambda^{r+1}-\lambda^{r}) \\ &= (\lambda^{r+1}-\lambda^{r})^T ((\Tilde{W} - W)\otimes I_{p}) (\lambda^{r+1}-\lambda^{r}) \\ &\geq (\lambda^{r+1}-\lambda^{r})^T \left(\frac{I_{np}+(\frac{1}{\rho}+1)W\otimes I_{p}}{\frac{1}{\rho}+2} - W\otimes I_{p} \right) (\lambda^{r+1}-\lambda^{r}) \\&= (\lambda^{r+1}-\lambda^{r})^T \left(\frac{I_{np}-W\otimes I_{p}}{\frac{1}{\rho}+2}  \right) (\lambda^{r+1}-\lambda^{r}) \\&\geq \frac{1-\lambda_{2}(W)}{\frac{1}{\rho}+2} ||\lambda^{r+1}-\lambda^{r}||^2 \\&= \frac{\rho (1-\lambda_{2}(W))}{1+2\rho} ||\lambda^{r+1}-\lambda^{r}||^2
\end{align*}
where the first inequality comes from Assumption \ref{AssM}. From \ref{dualdiff} we also have:
\begin{align*}
    & \ \ \ \ \ ||A^T (\lambda^{r+1}-\lambda^{r})||^2 \\
    &=||\nabla f(\mathbf{x}^{r}) - \nabla f(\mathbf{x}^{r-1}) +  C_{1}\mathbf{w}^{r+1} + \rho \left(I_{np}-\Tilde{W}\otimes I_{p}-((\Tilde{W} - W)\otimes I_{p})\right) (\mathbf{x}^{r} - \mathbf{x}^{r-1})||^2 \\
    &\le 2||\nabla f(\mathbf{x}^{r}) - \nabla f(\mathbf{x}^{r-1}) +\rho \left(I_{np}-\Tilde{W}\otimes I_{p}-((\Tilde{W} - W)\otimes I_{p})\right) (\mathbf{x}^{r} - \mathbf{x}^{r-1}) ||^2 + 2||C_{1}\mathbf{w}^{r+1}||^2 \\
    &\le 4||\nabla f(\mathbf{x}^{r}) - \nabla f(\mathbf{x}^{r-1})||^2 + 4\rho^2 ||\left(I_{np}-\Tilde{W}\otimes I_{p}-((\Tilde{W} - W)\otimes I_{p})\right) (\mathbf{x}^{r} - \mathbf{x}^{r-1})||^2 + 2||\mathbf{w}^{r+1}||_{C_{1}^T C_{1}}^2 \\
    &\le \left(4l^2 + 4\rho^2 ||I_{np}-\Tilde{W}\otimes I_{p}-((\Tilde{W} - W)\otimes I_{p})||^2 \right) ||\mathbf{x}^{r}-\mathbf{x}^{r-1}||^2 + 2||\mathbf{w}^{r+1}||_{C_{1}^T C_{1}}^2
\end{align*}
where the second and third inequality come from the Cauchy Schwartz inequality. Combine the above two inequalities together we have:
\begin{align*}
    &\frac{\rho (1-\lambda_{2}(W))}{1+2\rho} ||\lambda^{r+1}-\lambda^{r}||^2 \le \left(4l^2 + 4\rho^2 ||I_{np}-\Tilde{W}\otimes I_{p}-((\Tilde{W} - W)\otimes I_{p})||^2 \right) ||\mathbf{x}^{r}-\mathbf{x}^{r-1}||^2 + 2||\mathbf{w}^{r+1}||_{C_{1}^T C_{1}}^2 \\  &||\lambda^{r+1}-\lambda^{r}||^2 \le \frac{1+2\rho}{\rho (1-\lambda_{2}(W))}\left(4l^2 + 4\rho^2 ||I_{np}-\Tilde{W}\otimes I_{p}-((\Tilde{W} - W)\otimes I_{p})||^2 \right) ||\mathbf{x}^{r}-\mathbf{x}^{r-1}||^2 + \\& \ \ \ \ \ \ \ \ \ \ \ \ \ \ \ \ \ \ \ \ \ \frac{2(1+2\rho)}{\rho (1-\lambda_{2}(W))}||\mathbf{w}^{r+1}||_{C_{1}^T C_{1}}^2
\end{align*}
Multiply $\frac{1}{\rho}$ on both sides of this inequality we show \ref{lemma2}.
\hfill $\blacksquare$

\subsection{Proof of Lemma \ref{Le3}}

From \ref{MEXTRA}, $\forall r \geq 0$ and $\forall \mathbf{x}$ we have:
    \begin{align*}
        \langle \nabla f(\mathbf{x}^{r}) + A^T \lambda^{r} + \rho (I_{np}-\Tilde{W}\otimes I_{p}) \mathbf{x}^{r} + \beta (\mathbf{x}^{r+1}-\mathbf{x}^{r}), \mathbf{x}-\mathbf{x}^{r+1} \rangle \geq 0
    \end{align*}
    which implies
\begin{align*}
    \langle \nabla f(\mathbf{x}^{r}) + A^T \lambda^{r+1} + &\rho (I_{np}-\Tilde{W}\otimes I_{p})\mathbf{x}^r + \rho ((\Tilde{W} - W)\otimes I_{p}) (\mathbf{x}^{r}-\mathbf{x}^{r+1}) - \\ &\rho ((\Tilde{W} - W)\otimes I_{p}) \mathbf{x}^{r} + \beta (\mathbf{x}^{r+1}-\mathbf{x}^{r}), \mathbf{x}_{1}-\mathbf{x}^{r+1} \rangle \geq 0 \\
     \implies \langle \nabla f(\mathbf{x}^{r}) + A^T \lambda^{r+1} + &\rho \left(I_{np}-\Tilde{W}\otimes I_{p} -((\Tilde{W} - W)\otimes I_{p}) \right) \mathbf{x}^{r} + \\ &\rho ((\Tilde{W} - W)\otimes I_{p}) (\mathbf{x}^{r}-\mathbf{x}^{r+1}) + \beta (\mathbf{x}^{r+1}-\mathbf{x}^{r}), \mathbf{x}_{1}-\mathbf{x}^{r+1} \rangle \geq 0
\end{align*}
    where the first inequality comes from \ref{MEXTRA}. Since the above inequality is true for all $r > 0$ then we also have:
    \begin{align*}
        \langle \nabla f(\mathbf{x}^{r-1}) + A^T \lambda^{r} + \rho \left(I_{np}-\Tilde{W}\otimes I_{p} -((\Tilde{W} - W)\otimes I_{p}) \right) \mathbf{x}^{r-1} + &\rho ((\Tilde{W} - W)\otimes I_{p}) (\mathbf{x}^{r-1}-\mathbf{x}^{r}) + \\ &\beta (\mathbf{x}^{r}-\mathbf{x}^{r-1}), \mathbf{x}_{2}-\mathbf{x}^{r} \rangle \geq 0
    \end{align*}
    Let $\mathbf{x}_{1} = \mathbf{x}^{r}$ and $\mathbf{x}_{2}=\mathbf{x}^{r+1}$ and take the difference of the above two inequalities we have:
    \begin{align*}
        \langle &\nabla f(\mathbf{x}^{r-1})-\nabla f(\mathbf{x}^{r}) + A^T (\lambda^{r}-\lambda^{r+1}) + \rho \left(I_{np}-\Tilde{W}\otimes I_{p} -((\Tilde{W} - W)\otimes I_{p}) \right) (\mathbf{x}^{r-1}-\mathbf{x}^r) + \\
        &\rho ((\Tilde{W} - W)\otimes I_{p}) \{(\mathbf{x}^{r-1}-\mathbf{x}^{r})-(\mathbf{x}^{r}-\mathbf{x}^{r+1}) \} + \beta \{(\mathbf{x}^{r}-\mathbf{x}^{r-1})-(\mathbf{x}^{r+1}-\mathbf{x}^{r}) \},\mathbf{x}^{r+1}-\mathbf{x}^{r} \rangle \geq 0
    \end{align*}
    which implies:
    \begin{align}
        \langle A^T (\lambda^{r+1}-\lambda^{r}),\mathbf{x}^{r+1}-\mathbf{x}^{r} \rangle \le \langle &\nabla f(\mathbf{x}^{r-1})-\nabla f(\mathbf{x}^{r}) +
        \rho \left(I_{np}-\Tilde{W}\otimes I_{p} -((\Tilde{W} - W)\otimes I_{p}) \right) (\mathbf{x}^{r-1}-\mathbf{x}^r) -  \nonumber \\ &C_{1}\mathbf{w}^{r+1}, \mathbf{x}^{r+1}-\mathbf{x}^{r} \rangle \label{I1}
    \end{align}
    The LHS of \ref{I1} can be expressed as:
    \begin{align}
        \langle A^T (\lambda^{r+1}-\lambda^{r}),\mathbf{x}^{r+1}-\mathbf{x}^{r} \rangle &= \langle \lambda^{r+1}-\lambda^{r},A(\mathbf{x}^{r+1}-\mathbf{x}^{r}) \rangle \nonumber \\
        &= \rho \langle A\mathbf{x}^{r+1}, A\mathbf{x}^{r+1} - A\mathbf{x}^{r}\rangle \nonumber \\
        &= \rho ||A\mathbf{x}^{r+1}||^2 - \rho \langle A\mathbf{x}^{r+1}, A\mathbf{x}^{r}\rangle \nonumber \\
        &= \rho ||A\mathbf{x}^{r+1}||^2 - \frac{\rho}{2}\left(||A\mathbf{x}^{r+1}||^2 + ||A\mathbf{x}^{r}||^2 - ||A(\mathbf{x}^{r+1}-\mathbf{x}^{r})||^2 \right) \nonumber \\
        &= \frac{\rho}{2} ||A\mathbf{x}^{r+1}||^2 - \frac{\rho}{2} ||A\mathbf{x}^{r}||^2 + \frac{\rho}{2} ||A(\mathbf{x}^{r+1}-\mathbf{x}^{r})||^2 \nonumber \\
        &\geq \frac{\rho}{2} ||A\mathbf{x}^{r+1}||^2 - \frac{\rho}{2} ||A\mathbf{x}^{r}||^2 \label{L1LHS}
    \end{align}
    The RHS of \ref{I1} equals to:
    \begin{align}
        & \ \ \ \ \langle \nabla f(\mathbf{x}^{r-1})-\nabla f(\mathbf{x}^{r}) +
        \rho \left(I_{np}-\Tilde{W}\otimes I_{p} -((\Tilde{W} - W)\otimes I_{p}) \right) (\mathbf{x}^{r-1}-\mathbf{x}^r) -  \nonumber C_{1}\mathbf{w}^{r+1}, \mathbf{x}^{r+1}-\mathbf{x}^{r} \rangle \nonumber \\ &= \langle \nabla f(\mathbf{x}^{r-1})-\nabla f(\mathbf{x}^{r}), \mathbf{x}^{r+1}-\mathbf{x}^{r} \rangle + \langle \rho \left(I_{np}-\Tilde{W}\otimes I_{p} -((\Tilde{W} - W)\otimes I_{p}) \right) (\mathbf{x}^{r-1}-\mathbf{x}^r), \mathbf{x}^{r+1}-\mathbf{x}^{r} \rangle - \nonumber \\ & \ \ \ \ \ \ \ \ \ \ \ \ \ \ \ \ \ \ \ \ \ \ \ \ \ \ \ \ \ \ \ \ \ \ \ \ \ \ \ \ \ \ \ \ \ \ \ \ \ \ \ \ \ \ \ \ \ \ \ \ \ \ \ \ \ \ \ \ \ \ \ \ \ \ \ \ \ \ \ \ \ \ \ \ \ \ \ \ \ \ \ \ \ \ \ \ \ \ \ \ \ \  \langle C_{1}\mathbf{w}^{r+1}, \mathbf{x}^{r+1}-\mathbf{x}^{r} \rangle \nonumber \\
        &\le \frac{1}{2l} ||\nabla f(\mathbf{x}^{r-1})-\nabla f(\mathbf{x}^{r})||^2 + \frac{l}{2} ||\mathbf{x}^{r+1}-\mathbf{x}^{r}||^2 - \frac{1}{2} ||\mathbf{w}^{r+1}||_{C_{1}}^2 - \frac{1}{2} ||\mathbf{x}^{r+1}-\mathbf{x}^{r}||_{C_{1}}^2 + \nonumber \\ &\ \ \ \ \frac{1}{2} ||\mathbf{w}^{r+1}-(\mathbf{x}^{r+1}-\mathbf{x}^{r})||_{C_{1}}^2 + \langle \rho \left(I_{np}-\Tilde{W}\otimes I_{p} -((\Tilde{W} - W)\otimes I_{p}) \right) (\mathbf{x}^{r-1}-\mathbf{x}^r), \mathbf{x}^{r+1}-\mathbf{x}^{r} \rangle \nonumber \\ &\le \frac{l}{2} ||\mathbf{x}^{r}-\mathbf{x}^{r-1}||^2 + \frac{l}{2} ||\mathbf{x}^{r+1}-\mathbf{x}^{r}||^2 + \frac{1}{2} ||\mathbf{x}^{r}-\mathbf{x}^{r-1}||_{C_{1}}^2 - \frac{1}{2} ||\mathbf{w}^{r+1}||_{C_{1}}^2 - \frac{1}{2} ||\mathbf{x}^{r+1}-\mathbf{x}^{r}||_{C_{1}}^2 + \nonumber \\ & \ \ \ \ \langle \rho \left(I_{np}-\Tilde{W}\otimes I_{p} -((\Tilde{W} - W)\otimes I_{p}) \right) (\mathbf{x}^{r-1}-\mathbf{x}^r), \mathbf{x}^{r+1}-\mathbf{x}^{r} \rangle \label{I2}
    \end{align}
    Recall $C_{2} = \rho \left(I_{np}-\Tilde{W}\otimes I_{p} -((\Tilde{W} - W)\otimes I_{p}) \right)$ then the last term of \ref{I2} equals to:
    \begin{align}
        & \ \ \ \ \langle \rho \left(I_{np}-\Tilde{W}\otimes I_{p} -((\Tilde{W} - W)\otimes I_{p}) \right) (\mathbf{x}^{r-1}-\mathbf{x}^r), \mathbf{x}^{r+1}-\mathbf{x}^{r} \rangle \nonumber \\
        &= \langle C_{2} (\mathbf{x}^{r-1}-\mathbf{x}^r), (\mathbf{x}^{r+1}-\mathbf{x}^{r}) + (\mathbf{x}^{r-1}-\mathbf{x}^{r}) - (\mathbf{x}^{r-1}-\mathbf{x}^{r}) \rangle \nonumber \\
        &= \langle C_{2} (\mathbf{x}^{r-1}-\mathbf{x}^{r}), \mathbf{w}^{r+1} \rangle - ||\mathbf{x}^{r-1}-\mathbf{x}^{r}||_{C_{2}}^2 \nonumber \\
        &\le \langle C_{2} (\mathbf{x}^{r-1}-\mathbf{x}^{r}), \mathbf{w}^{r+1} \rangle \nonumber \\
        &= -\frac{1}{2} ||\mathbf{w}^{r+1}-(\mathbf{x}^{r-1}-\mathbf{x}^{r})||_{C_{2}}^2 + \frac{1}{2} ||\mathbf{x}^{r}-\mathbf{x}^{r-1}||_{C_{2}}^2 + \frac{1}{2} ||\mathbf{w}^{r+1}||_{C_{2}}^2 \nonumber \\  &= -\frac{1}{2}||\mathbf{x}^{r+1}-\mathbf{x}^{r}||_{C_{2}}^2 + \frac{1}{2} ||\mathbf{x}^{r}-\mathbf{x}^{r-1}||_{C_{2}}^2 + \frac{1}{2} ||\mathbf{w}^{r+1}||_{C_{2}}^2 \label{I3}
    \end{align}
    Then combine the above four inequalities \ref{I1},\ref{L1LHS},\ref{I2},\ref{I3} together we have:
    \begin{align*}
        \frac{\rho}{2} ||A\mathbf{x}^{r+1}||^2 - \frac{\rho}{2} ||A\mathbf{x}^{r}||^2 &\le \frac{l}{2} ||\mathbf{x}^{r}-\mathbf{x}^{r-1}||^2 + \frac{l}{2} ||\mathbf{x}^{r+1}-\mathbf{x}^{r}||^2 + \frac{1}{2} ||\mathbf{x}^{r}-\mathbf{x}^{r-1}||_{C_{1}}^2 - \frac{1}{2} ||\mathbf{w}^{r+1}||_{C_{1}}^2 \\ &- \frac{1}{2} ||\mathbf{x}^{r+1}-\mathbf{x}^{r}||_{C_{1}}^2 - \frac{1}{2} ||\mathbf{x}^{r+1}-\mathbf{x}^{r}||_{C_{2}}^2 + \frac{1}{2} ||\mathbf{x}^{r}-\mathbf{x}^{r-1}||_{C_{2}}^2 + \frac{1}{2} ||\mathbf{w}^{r+1}||_{C_{2}}^2
    \end{align*}
    which implies:
    \begin{align*}
        &\ \ \ \ \frac{\rho}{2} ||A\mathbf{x}^{r+1}||^2 + \frac{1}{2} ||\mathbf{x}^{r+1}-\mathbf{x}^{r}||_{C_{1}}^2 + \frac{1}{2} ||\mathbf{x}^{r+1}-\mathbf{x}^{r}||_{C_{2}}^2 \\ & \le \frac{\rho}{2} ||A\mathbf{x}^{r}||^2 + \frac{1}{2} ||\mathbf{x}^{r}-\mathbf{x}^{r-1}||_{C_{1}}^2 + \frac{1}{2} ||\mathbf{x}^{r}-\mathbf{x}^{r-1}||_{C_{2}}^2 - \frac{1}{2} ||\mathbf{w}^{r+1}||_{C-C_{2}}^2 + \\ & \ \ \ \ \frac{l}{2} ||\mathbf{x}^{r}-\mathbf{x}^{r-1}||^2 + \frac{l}{2} ||\mathbf{x}^{r+1}-\mathbf{x}^{r}||^2
    \end{align*}
    Recall that $C_{1} = \beta I - \rho ((\Tilde{W} - W)\otimes I_{p})$, $C_{2} = \rho \left(I_{np}-\Tilde{W}\otimes I_{p} -((\Tilde{W} - W)\otimes I_{p}) \right)$ and $C_{3}=\beta I_{np} - \rho (I_{np}-\Tilde{W}\otimes I_{p})$ then we have:
    \begin{align*}
        C_{1} - C_{2} &= \beta I_{np} - \rho ((\Tilde{W} - W)\otimes I_{p}) - \rho \left(I_{np}-\Tilde{W}\otimes I_{p} -((\Tilde{W} - W)\otimes I_{p}) \right) \\ &= \beta I_{np} - \rho ((\Tilde{W} - W)\otimes I_{p}) - \rho (I_{np}-\Tilde{W}\otimes I_{p}) + \rho ((\Tilde{W} - W)\otimes I_{p}) \\ &= \beta I_{np} - \rho (I_{np}-\Tilde{W}\otimes I_{p}) \\ &=C_{3}
    \end{align*}
    Therefore, we have:
    \begin{align*}
         &\ \ \ \ \frac{\rho}{2} ||A\mathbf{x}^{r+1}||^2 + \frac{1}{2} ||\mathbf{x}^{r+1}-\mathbf{x}^{r}||_{C_{1}}^2 + \frac{1}{2} ||\mathbf{x}^{r+1}-\mathbf{x}^{r}||_{C_{2}}^2 \\ & \le \frac{\rho}{2} ||A\mathbf{x}^{r}||^2 + \frac{1}{2} ||\mathbf{x}^{r}-\mathbf{x}^{r-1}||_{C_{1}}^2 + \frac{1}{2} ||\mathbf{x}^{r}-\mathbf{x}^{r-1}||_{C_{2}}^2 - \frac{1}{2} ||\mathbf{w}^{r+1}||_{C_{3}}^2 + \\ & \ \ \ \ \frac{l}{2} ||\mathbf{x}^{r}-\mathbf{x}^{r-1}||^2 + \frac{l}{2} ||\mathbf{x}^{r+1}-\mathbf{x}^{r}||^2
    \end{align*} \hfill $\blacksquare$

\subsection{Proof of Lemma \ref{Le4}}
\begin{proof}
    Recall the Lyapunov function:
\begin{align}
    P_{c} (\mathbf{x}^{r+1},\mathbf{x}^{r},\lambda^{r+1}) &= L_{\rho}(\mathbf{x}^{r+1},\lambda^{r+1}) + \frac{c\rho}{2} \|A\mathbf{x}^{r+1}\|^2 + \frac{c}{2} \|\mathbf{x}^{r+1}-\mathbf{x}^{r}\|^2_{C_{1}} + \frac{c}{2} \|\mathbf{x}^{r+1}-\mathbf{x}^{r}\|^2_{C_{2}} + \nonumber \\ &\ \ \ \left\{\frac{cl}{2} + \frac{1+2\rho}{\rho^2 (1-\lambda_{2}(W))} \left(4l^2 + 4\rho^2 \|(I_{n}\!-\!2\Tilde{W}\!+W)\otimes I_{p}\|^2 \right)\right\} \|\mathbf{x}^{r+1}-\mathbf{x}^{r}\|^2  \label{lyafun}
\end{align} where $C_{1} = \beta I_{np} - \rho A^T A$, $ C_{2} = \rho (I_{n}-\Tilde{W} -(\Tilde{W}-W) )\otimes I_{p}$ and $A = ((\Tilde{W} - W)\otimes I_{p})^{1/2}$. In this Lemma, we will show that $P_{c}$ is monotonically decrease. Combining Lemmas~\ref{Le1},~\ref{Le2}, and~\ref{Le3} while multiplying both sides of the inequality derived in Lemma~\ref{Le3} by constant $c > 0$, we get:
\begin{align*}
    & \ \ \ \ \ \ L_{\rho}(\mathbf{x}^{r+1},\lambda^{r+1}) + \frac{c\rho}{2} \|A\mathbf{x}^{r+1}\|^2 + \frac{c}{2} \|\mathbf{x}^{r+1}-\mathbf{x}^{r}\|^2_{C_{1}} + \frac{c}{2} \|\mathbf{x}^{r+1}-\mathbf{x}^{r}\|^2_{C_{2}} \\
    &= L_{\rho}(\mathbf{x}^{r},\lambda^{r}) + \frac{c\rho}{2} \|A\mathbf{x}^{r}\|^2 + \frac{c}{2} \|\mathbf{x}^{r}-\mathbf{x}^{r-1}\|^2_{C_{1}} + \frac{c}{2} \|\mathbf{x}^{r}-\mathbf{x}^{r-1}\|^2_{C_{2}} + \\ &\ \ \ \ 
   \left\{\frac{cl}{2} + \frac{1+2\rho}{\rho^2 (1-\lambda_{2}(W))} \left(4l^2 + 4\rho^2 \|(I_{n}\!-\!2\Tilde{W}\!+W)\otimes I_{p}\|^2 \right)\right\} \|\mathbf{x}^{r}-\mathbf{x}^{r-1}\|^2 - \\ & \ \ \ \ \|\mathbf{w}^{r+1}\|_{\frac{c}{2}C_{3} - \frac{2(1+2\rho)}{\rho^2 (1-\lambda_{2}(W))} C_{1}^T C_{1}}^2 - \left(\beta - \frac{L}{2} -\frac{cl}{2}\right) \|\mathbf{x}^{r+1}-\mathbf{x}^{r}\|^2.
\end{align*}
Hence,
\begin{align*}
    P_{c} (\mathbf{x}^{r+1},\mathbf{x}^{r},&\lambda^{r+1}) \le P_{c} (\mathbf{x}^{r},\mathbf{x}^{r-1},\lambda^{r}) - \|\mathbf{w}^{r+1}\|_{\frac{c}{2}C_{3} - \frac{2(1+2\rho)}{\rho^2 (1-\lambda_{2}(W))}C_{1}^T C_{1}}^2 - \\ & \left\{ \beta - cl - \frac{L}{2} - \frac{1+2\rho}{\rho^2 (1-\lambda_{2}(W))} \left(4l^2 + 4\rho^2 \|I_{n}\!-\!2\Tilde{W}\!+W)\otimes I_{p}\|^2 \right) \right\} \|\mathbf{x}^{r+1}-\mathbf{x}^{r}\|^2
\end{align*}

    Let $\lambda$ be an arbitrary eigenvalue of $(\Tilde{W}-W)\otimes I_{p}$ and from Assumption~\ref{AssM} we know $\lambda \geq 0$. 
    Since $\rho > \max\left\{ 1 + \lambda_{\max} (\Tilde{W}-W), \frac{8l + \sqrt{64l^2 + 16l(1-\lambda_{2}(W))}}{2(1-\lambda_{2}(W))} \right\}$ then we have $\rho^2 (1-\lambda_{2}(W)) > 4l(1+2\rho)$ which implies $\frac{\rho^2 (1-\lambda_{2}(W))}{4l(1+2\rho)} > 1$. Therefore, we can always find a constant $a$ such that $1 < a \le \frac{\rho^2 (1-\lambda_{2}(W))}{4l(1+2\rho)}$ which is equivalent to $l < al \le \frac{\rho^2 (1-\lambda_{2}(W))}{4(1+2\rho)}$. Let $c = \frac{\beta}{al}$ then we have $c = \frac{\beta}{al} \geq \frac{4(1+2\rho)}{\rho^2 (1-\lambda_{2}(W))} \beta$ which implies:
    \begin{align*}
        c &> \frac{4(1+2\rho)}{\rho^2 (1-\lambda_{2}(W))}\beta - \frac{4(1+2\rho)}{\rho^2 (1-\lambda_{2}(W))} \left((\rho + 1)\lambda - 2\lambda - \lambda^2 \right) \\
        &= \frac{4(1+2\rho)}{\rho^2 (1-\lambda_{2}(W))} \left(\beta -  (\rho + 1)\lambda + 2\lambda + \lambda^2\right) \\
        &\geq \frac{4(1+2\rho)}{\rho^2 (1-\lambda_{2}(W))} \left(\beta -  (\rho + 1)\lambda + 2\lambda + \frac{\lambda^2}{\beta - \rho \lambda -\lambda} \right) \\
        &= \frac{4(1+2\rho)}{\rho^2 (1-\lambda_{2}(W))} \left\{\frac{(\beta -  (\rho + 1)\lambda)^2}{\beta - \rho \lambda -\lambda} + \frac{2\lambda (\beta - \rho \lambda -\lambda)}{\beta - \rho \lambda -\lambda} + \frac{\lambda^2}{\beta - \rho \lambda -\lambda} \right\} \\
        &= \frac{4(1+2\rho)}{\rho^2 (1-\lambda_{2}(W))} \frac{(\beta - \rho \lambda)^2}{\beta - \rho \lambda -\lambda}
    \end{align*}
    where the first inequality comes from $\rho > 1 + \lambda_{\max}(\Tilde{W}-W) \geq 2 + \lambda -1$ which implies $(\rho+1)\lambda > 2\lambda +\lambda^2$ and the second inequality comes from $\beta - \rho \lambda - \lambda \geq 1$. Since $c > \frac{4(1+2\rho)}{\rho^2 (1-\lambda_{2}(W))} \frac{(\beta - \rho \lambda)^2}{\beta - \rho \lambda -\lambda}$ then we have:
    \begin{align}
        &\frac{c}{2} (\beta - \rho \lambda - \lambda) > \frac{2(1+2\rho)}{\rho^2 (1-\lambda_{2}(W))} (\beta - \rho \lambda)^2 \nonumber \\
        \implies &\frac{c}{2} (\beta - \rho \lambda - \lambda) - \frac{2(1+2\rho)}{\rho^2 (1-\lambda_{2}(W))} (\beta - \rho \lambda)^2 > 0 \label{eigen}
    \end{align}
    Recall that $C_{1} = \beta I_{np} - \rho (\Tilde{W}-W)\otimes I_{p}$ and $\lambda$ is an arbitrary eigenvalue of $(\Tilde{W}-W)\otimes I_{p}$ then the LHS of \eqref{eigen} is an eigenvalue of $\frac{c}{2} (\beta I_{np} - \rho (\Tilde{W}-W)\otimes I_{p} - (\Tilde{W}-W))\otimes I_{p} - \frac{2(1+2\rho)}{\rho^2 (1-\lambda_{2}(W))} C_{1}^T C_{1}$. Since \eqref{eigen} is true for all $\lambda$ then we have:
    \begin{align*}
        \frac{c}{2} \left\{\beta I_{np} - \rho (\Tilde{W}-W)\otimes I_{p} - (\Tilde{W}-W)\otimes I_{p}\right\} - \frac{2(1+2\rho)}{\rho^2 (1-\lambda_{2}(W))} C_{1}^T C_{1} \succ 0
    \end{align*}
    From Assumption \ref{AssM} we have $I_{n} - \Tilde{W} \preceq (\frac{1}{\rho} + 1)(\Tilde{W}-W)$ which implies:
    \begin{align*}
        &\ \ \ \ \ \frac{c}{2} C_{3} - \frac{2(1+2\rho)}{\rho^2 (1-\lambda_{2}(W))}C_{1}^T C_{1} \\
        & =\frac{c}{2}\left( \beta I_{np} - \rho (I-\Tilde{W})\otimes I_{p}\right) - \frac{2(1+2\rho)}{\rho^2 (1-\lambda_{2}(W))}C_{1}^T C_{1} \\ &\succeq \frac{c}{2}\left( \beta I - \rho (\frac{1}{\rho} + 1)(\Tilde{W}-W)\otimes I_{p}\right) - \frac{2(1+2\rho)}{\rho^2 (1-\lambda_{2}(W))}C_{1}^T C_{1} \\
        &= \frac{c}{2} \left\{\beta I - \rho (\Tilde{W}-W)\otimes I_{p} - (\Tilde{W}-W)\otimes I_{p}\right\} - \frac{2(1+2\rho)}{\rho^2 (1-\lambda_{2}(W))} C_{1}^T C_{1} \\ 
        &\succ 0
    \end{align*}
    Since $a > 1$ then $1-\frac{1}{a} > 0$ therefore if $\beta$ satisfies:
    \begin{align*}
        \beta > \left\{\frac{L}{2} + \frac{1+2\rho}{\rho^2 (1-\lambda_{2}(W))} \left(4l^2 + 4\rho^2 \|I-\Tilde{W}-(\Tilde{W}-W)\|^2 \right)\right\} \Biggm/ (1-1/a)
    \end{align*}
    then we have:
    \begin{align*}
        \beta - cl - \frac{L}{2} - \frac{1+2\rho}{\rho^2 (1-\lambda_{2}(W))} \left(4l^2 + 4\rho^2 \|I-\Tilde{W}-(\Tilde{W}-W)\|^2 \right) > 0
    \end{align*}
\end{proof}

\subsection{Proof of Lemma \ref{Le5}}
\begin{proof}
From \eqref{primal-obj} we have:
    \begin{align*}
        L_{\rho}(\mathbf{x}^{r+1},\lambda^{r+1}) &= f(\mathbf{x}^{r+1}) + \langle \lambda^{r+1},A\mathbf{x}^{r+1} \rangle + \frac{\rho}{2} \|\mathbf{x}^{r+1}\|_{I-\Tilde{W}}^2 \\
        &= f(\mathbf{x}^{r+1}) + \frac{1}{\rho} \langle \lambda^{r+1},\lambda^{r+1}-\lambda^{r} \rangle + \frac{\rho}{2} \|\mathbf{x}^{r+1}\|_{I-\Tilde{W}}^2 \\
        &= f(\mathbf{x}^{r+1}) + \frac{1}{2\rho} \left(\|\lambda^{r+1}\|^2 - \|\lambda^{r}\|^2 + \|\lambda^{r+1}-\lambda^{r}\|^2 \right) + \frac{\rho}{2} \|\mathbf{x}^{r+1}\|_{I-\Tilde{W}}^2 \\
        which \ implies \ that \\
         \sum_{r=1}^{T} L_{\rho}(\mathbf{x}^{r+1},\lambda^{r+1}) &=  \sum_{r=1}^{T} \left\{f(\mathbf{x}^{r+1}) + \frac{1}{2\rho} \|\lambda^{r+1}-\lambda^{r}\|^2 + \frac{\rho}{2} \|\mathbf{x}^{r+1}\|_{I-\Tilde{W}}^2 +\frac{\rho}{2}\|\lambda^{r+1}\|^2 - \frac{\rho}{2}\|\lambda^{r}\|^2\right\} \\
        &= \sum_{r=1}^{T} \left\{f(\mathbf{x}^{r+1}) + \frac{1}{2\rho} \|\lambda^{r+1}-\lambda^{r}\|^2 + \frac{\rho}{2} \|\mathbf{x}^{r+1}\|_{I-\Tilde{W}}^2\right\} + \frac{\rho}{2}\|\lambda^{T+1}\|^2 - \frac{\rho}{2}\|\lambda^{1}\|^2 \\ &\geq - \frac{\rho}{2}\|\lambda^{1}\|^2
    \end{align*}
    From the definition of $P_{c}$ in (\ref{lyafun}), we can see that except for the $L_{\rho}$ part, the rest terms are all non-negative, therefore we have:
    \begin{align*}
        \sum_{r=1}^{T} P_{c} (\mathbf{x}^{r+1},\mathbf{x}^{r},\lambda^{r+1}) \geq - \frac{\rho}{2}\|\lambda^{1}\|^2
    \end{align*}
    Moreover, if $\beta$, $\rho$ satisfy \eqref{conr1}, \eqref{conr2}, \eqref{conb1} and \eqref{conb2} then $P_{c}$ is monotonically decrease and therefore upper bounded by $P_{c} (\mathbf{x}^{1},\mathbf{x}^{0},\lambda^{1})$. Therefore, $P_{c} (\mathbf{x}^{r+1},\mathbf{x}^{r},\lambda^{r+1})$ is uniformly lower bounded which completes the proof.
\end{proof}
\subsection{Proof of Theorem \ref{Thm}}
\begin{proof}
    
Form Lemmas \ref{lemma4} and  \ref{lemma5} we know $P_{c} (\mathbf{x}^{r+1},\mathbf{x}^{r},\lambda^{r+1})$ is monotonically decrease and lower bounded. Therefore, the limit of $P_{c} (\mathbf{x}^{r+1},\mathbf{x}^{r},\lambda^{r+1})$ exists which implies that the difference between two consecutive elements approaches 0 i.e. $\lim_{r\rightarrow \infty} \mathbf{w}^{r+1} = 0$ and $\lim_{r\rightarrow \infty} \mathbf{x}^{r+1}-\mathbf{x}^{r} = 0$. Therefore, we have $\lim_{r\rightarrow \infty} \|\lambda^{r+1}-\lambda^{r}\|^2 = 0$ (see \eqref{lemma2}) which directly implies $\lim_{r\rightarrow \infty} \|A\mathbf{x}^{r+1}\|^2 = 0$. Since $\textup{Null}\{A \} = \textup{span}\{\bm{1}\}$ then we have the asymptotically consensus i.e. $\lim_{r\rightarrow \infty} \|\mathbf{x}^{r}-\Bar{\mathbf{x}}^{r}\|^2 = 0$.
    From \eqref{MEXTRA} we have:
    \begin{align*}
        &\nabla f(\mathbf{x}^{r}) + A^T \lambda^{r} + \rho (I-\Tilde{W}) \mathbf{x}^{r} + \beta (\mathbf{x}^{r+1}-\mathbf{x}^{r}) = 0 \\
        \implies & \bm 1^T (\nabla f(\mathbf{x}^{r}) + A^T \lambda^{r} + \rho (I-\Tilde{W}) \mathbf{x}^{r} + \beta (\mathbf{x}^{r+1}-\mathbf{x}^{r})) = 0 \\
        \implies &\bm 1^T (\nabla f(\mathbf{x}^{r}) + \beta (\mathbf{x}^{r+1}-\mathbf{x}^{r})) = 0 \\
        \implies & \lim_{r\rightarrow \infty} \bm 1^T (\nabla f(\mathbf{x}^{r}) + \beta (\mathbf{x}^{r+1}-\mathbf{x}^{r})) = 0 \\
        \implies & \lim_{r\rightarrow \infty} \bm 1^T \nabla f(\mathbf{x}^{r}) = 0
    \end{align*}
    where the third equality comes from Assumption \ref{AssM} and the last equality comes from the fact that $\lim_{r\rightarrow \infty} \mathbf{x}^{r                  +1}-\mathbf{x}^{r} = 0$.

    To show the sub-linear convergence rate, we recall $T_{\epsilon} = \min\{r: \|\nabla f(\mathbf{x}^{r})+A^T \lambda^{r}+((I_{n}-\Tilde{W})\otimes I_{p})\mathbf{x}^{r}\|^2 + \|A\mathbf{x}^{r}\|^2 \le \epsilon \}$.
    From Lemma \ref{Le2}, there exists two positive constants $a_{1},b_{1}$ such that:
    \begin{align*}
        ||A\mathbf{x}^{r+1}||^2 \le a_{1}||\mathbf{x}^{r}-\mathbf{x}^{r-1}||^2 + b_{1} ||(\mathbf{x}^{r+1}-\mathbf{x}^{r})-(\mathbf{x}^{r}-\mathbf{x}^{r-1})||^2
    \end{align*}
    which implies:
    \begin{align*}
        ||A\mathbf{x}^{r+1}|| &\le \sqrt{a_{1}}||\mathbf{x}^{r}-\mathbf{x}^{r-1}|| + \sqrt{b_{1}}||(\mathbf{x}^{r+1}-\mathbf{x}^{r})-(\mathbf{x}^{r}-\mathbf{x}^{r-1})|| \\
        & = \sqrt{a_{1}}||\mathbf{x}^{r}-\mathbf{x}^{r-1}|| - \sqrt{a_{1}}||\mathbf{x}^{r+1}-\mathbf{x}^{r}|| + \sqrt{a_{1}}||\mathbf{x}^{r+1}-\mathbf{x}^{r}|| + \sqrt{b_{1}}||(\mathbf{x}^{r+1}-\mathbf{x}^{r})-(\mathbf{x}^{r}-\mathbf{x}^{r-1})|| \\
        & \le \sqrt{a_{1}}||(\mathbf{x}^{r+1}-\mathbf{x}^{r})-(\mathbf{x}^{r}-\mathbf{x}^{r-1})|| + \sqrt{a_{1}}||\mathbf{x}^{r+1}-\mathbf{x}^{r}|| + \sqrt{b_{1}}||(\mathbf{x}^{r+1}-\mathbf{x}^{r})-(\mathbf{x}^{r}-\mathbf{x}^{r-1})|| \\
        & = (\sqrt{a_{1}}+\sqrt{b_{1}})||(\mathbf{x}^{r+1}-\mathbf{x}^{r})-(\mathbf{x}^{r}-\mathbf{x}^{r-1})|| + \sqrt{a_{1}}||\mathbf{x}^{r+1}-\mathbf{x}^{r}||
    \end{align*}
    which implies:
    \begin{align*}
        ||A\mathbf{x}^{r+1}||^2 &\le 2(\sqrt{a_{1}}+\sqrt{b_{1}})^2||(\mathbf{x}^{r+1}-\mathbf{x}^{r})-(\mathbf{x}^{r}-\mathbf{x}^{r-1})||^2 + 2a_{1}||\mathbf{x}^{r+1}-\mathbf{x}^{r}||^2 \\
        & \le max\{2(\sqrt{a_{1}}+\sqrt{b_{1}})^2, 2a_{1}\} (||\mathbf{w}^{r+1}||^2 + ||\mathbf{x}^{r+1}-\mathbf{x}^{r}||^2) \\
        & := c_{1} (||\mathbf{w}^{r+1}||^2 + ||\mathbf{x}^{r+1}-\mathbf{x}^{r}||^2)
    \end{align*}
    where $\mathbf{w}^{r+1} = (\mathbf{x}^{r+1}-\mathbf{x}^{r}) - (\mathbf{x}^{r}- \mathbf{x}^{r-1})$. Since:
    \begin{align*}
        ||A\mathbf{x}^{r}|| \le ||A\mathbf{x}^{r+1}-A\mathbf{x}^{r}|| + ||A\mathbf{x}^{r+1}||
    \end{align*}
    then 
    \begin{align*}
        ||A\mathbf{x}^{r}||^2 &\le 2||A\mathbf{x}^{r+1}-A\mathbf{x}^{r}||^2 + 2||A\mathbf{x}^{r+1}||^2 \\
        ||A\mathbf{x}^{r}||^2 &\le 2||\mathbf{x}^{r+1}-\mathbf{x}^{r}||^2 + 2||A\mathbf{x}^{r+1}||^2 \\
        & \le (2+2c_{1})||\mathbf{x}^{r+1}-\mathbf{x}^{r}||^2 + 2c_{1}||\mathbf{w}^{r+1}||^2 \\
        & \le (2+2c_{1})(||\mathbf{x}^{r+1}-\mathbf{x}^{r}||^2 + ||\mathbf{w}^{r+1}||^2) \\
        & := d_{1}(||\mathbf{x}^{r+1}-\mathbf{x}^{r}||^2 + ||\mathbf{w}^{r+1}||^2) 
    \end{align*}
    where the second inequality comes from $||A||\le 1$. From $\mathbf{x}$ iterates updates in \eqref{MEXTRA} we have:
    \begin{align*}
        ||\nabla f(\mathbf{x}^{r}) + A^T \lambda^{r} + \rho (I-\Tilde{W}) \mathbf{x}^{r}||^2 = \beta ||\mathbf{x}^{r+1}-\mathbf{x}^{r}||^2
    \end{align*}
    From lemma \ref{lemma4}, there exist $e_{1},e_{2} > 0$ such that:
    \begin{align*}
        e_{1}||\mathbf{x}^{r+1}-\mathbf{x}^{r}||^2 + e_{2}||\mathbf{w}^{r+1}||^2 \le P_{c}(\mathbf{x}^{r},\mathbf{x}^{r-1},\lambda^{r})- P_{c}(\mathbf{x}^{r+1},\mathbf{x}^{r},\lambda^{r+1}) 
    \end{align*}
    which implies
    \begin{align*}
        min\{e_{1},e_{2} \} (||\mathbf{x}^{r+1}-\mathbf{x}^{r}||^2 + ||\mathbf{w}^{r+1}||^2) \le P_{c}(\mathbf{x}^{r},\mathbf{x}^{r-1},\lambda^{r})- P_{c}(\mathbf{x}^{r+1},\mathbf{x}^{r},\lambda^{r+1})  
    \end{align*}
    Therefore, we have:
    \begin{align*}
        ||A\mathbf{x}^{r}||^2 + ||\nabla f(\mathbf{x}^{r}) + A^T \lambda^{r} + \rho (I-\Tilde{W}) \mathbf{x}^{r}||^2 \le \frac{d_{1}+\beta}{min\{e_{1}, e_{2} \}} (P_{c}(\mathbf{x}^{r},\mathbf{x}^{r-1},\lambda^{r})- P_{c}(\mathbf{x}^{r+1},\mathbf{x}^{r},\lambda^{r+1}))
    \end{align*}
    By definition of $T_{\epsilon}$ and $P_{c}(\mathbf{x}^{r+1},\mathbf{x}^{r},\lambda^{r+1}) >$ \underline{$p$} (see Lemma \ref{Le5}), we have:
    \begin{align*}
        T_{\epsilon}\epsilon &\le \frac{d_{1}+\beta}{min\{e_{1},e_{2} \}} \sum_{r=1}^{T} (P_{c}(\mathbf{x}^{r},\mathbf{x}^{r-1},\lambda^{r})- P_{c}(\mathbf{x}^{r+1},\mathbf{x}^{r},\lambda^{r+1})) \\
        & = \frac{d_{1}+\beta}{min\{e_{1}, e_{2} \}} (P_{c}(\mathbf{x}^{1},\mathbf{x}^{0},\lambda^{1})- P_{c}(\mathbf{x}^{T+1},\mathbf{x}^{T},\lambda^{T+1})) \\
        &\le \frac{d_{1}+\beta}{min\{e_{1}, e_{2} \}} (P_{c}(\mathbf{x}^{1},\mathbf{x}^{0},\lambda^{1})- \underline{p})
    \end{align*}
    Therefore, we have $T_{\epsilon} = \mathcal{O}(\frac{1}{\epsilon})$.
\end{proof}
\subsection{Proof of Lemma \ref{pro}}
\begin{proof}
    \begin{align}
        & \ \ \ \ \ \frac{I+(\frac{1}{\rho_{1}}+1)W}{\frac{1}{\rho_{1}}+2} - \frac{I + (\frac{1}{\rho_{2}}+1)W}{\frac{1}{\rho_{2}}+2}\nonumber \\ &= \frac{\rho_{1}I+(1+\rho_{1})W}{1+2\rho_{1}} - \frac{\rho_{2}I+(1+\rho_{2})W}{1+2\rho_{2}} \nonumber \\
        &= \frac{(1+2\rho_{2})\rho_{1}I+(1+\rho_{1})(1+2\rho_{2})W-(1+2\rho_{1})\rho_{2}I-(1+\rho_{2})(1+2\rho_{1})W}{(1+2\rho_{1})(1+2\rho_{2})}\nonumber \\
        &= \frac{(\rho_{1}-\rho_{2})I + (\rho_{2}-\rho_{1})W}{(1+2\rho_{1})(1+2\rho_{2})} .\label{MM}
    \end{align}
    Let $\lambda$ to be an arbitrary eigenvalue of $W$, from assumption \ref{AssM} we know that $\lambda \le 1$. Then the eigenvalue of \eqref{MM} becomes:
    \begin{align*}
        \frac{\rho_{1}-\rho_{2} + (\rho_{2}-\rho_{1})\lambda}{(1+2\rho_{1})(1+2\rho_{2})} = \frac{(\rho_{1}-\rho_{2})(1-\lambda)}{(1+2\rho_{1})(1+2\rho_{2})} \geq 0,
    \end{align*}
    where the inequality comes from $\rho_{1} \geq \rho_{2}$ and $\lambda \le 1$.
\end{proof}
\subsection{Proof of Lemma \ref{prop3}}
\begin{proof}
Since $\Tilde{W} = aI + bW$, then we have:
    \begin{align*}
        \Tilde{W}-W &= aI + bW - W \\&= aI - aW + (a+b-1)W \\
        &= a(I-W).
    \end{align*}
    Therefore, $\textup{Null}(\Tilde{W}-W) = \textup{Null}(I-W) = \textup{span} \{\bm 1 \}$.
    \begin{align*}
        (I-\Tilde{W}) \bm 1 &= (I - aI - bW) \bm 1 \\
        &=\bm 1 - a\bm 1 -b\bm 1\\
        &=0,
    \end{align*}
    which completes the proof.
\end{proof}
\end{document}